\documentclass[oneside,12pt,leqno]{amsart}
\usepackage{amssymb}
\swapnumbers
\theoremstyle{plain}
\newtheorem{thm}{Theorem}[section]
\newtheorem{prop}[thm]{Proposition}
\newtheorem{cor}[thm]{Corollary}
\newtheorem{lemma}[thm]{Lemma}
\newtheorem{maintheo}[thm]{Main Theorem}
\theoremstyle{definition}
\newtheorem{defn}[thm]{Definition}
\newtheorem{anz}[thm]{Ansatz}
\newtheorem{notn}[thm]{Notation}
\theoremstyle{remark}
\newtheorem{remk}[thm]{Remark}
\newcommand{\proofof}[1]{\end{#1}\begin{proof}}
\newcommand{\emphdef}{\textit}
\numberwithin{equation}{section}

\newcommand{\thismonth}{\ifcase\month\or
  January\or February\or March\or April\or May\or June\or July\or
  August\or September\or October\or November\or December\fi
  \space\number\year}
\newcommand{\tcite}[1]{\textup{\cite{#1}}}
\DeclareMathAlphabet{\mathrmsl}{OT1}{cmr}{m}{sl}
\newcommand{\symb}[2]{\newcommand{#1}{\mathit{#2}} }
\newcommand{\rssymb}[2]{\newcommand{#1}{\mathrmsl{#2}} }
\newcommand{\oper}[3][n]{\newcommand{#2}{\mathop{\mathrm{#3}}%
\ifx n#1\nolimits\else\limits\fi} }
\newcommand{\rsoper}[3][n]{\newcommand{#2}{\mathop{\mathrmsl{#3}}%
\ifx n#1\nolimits\else\limits\fi} }
\oper\image{im}
\oper\End{End}
\oper\Aut{Aut}
\oper\Ad{Ad}
\oper\Alt{Alt}
\oper\Sym{Sym}
\oper[l]\Res{Res}
\rsoper{\sym}{sym}
\rsoper{\alt}{alt}
\rsoper\trace{tr}
\rsoper{\ptr}{ptr}
\rsoper{\divg}{div}
\symb\ev{ev}
\rssymb\iden{id}
\rssymb\Scal{scal}
\newcommand{\ip}[1]{\langle#1\rangle}

\newcommand{\lie}[1]{\mathfrak{#1}}
\newcommand{\N}{{\mathbb N}}

\newcommand{\R}{{\mathbb R}}

\newcommand{\eps}{\varepsilon}
\renewcommand{\geq}{\geqslant}
\renewcommand{\leq}{\leqslant}
\newcommand{\Dsum}{\bigoplus}
\newcommand{\dsum}{\oplus}
\newcommand{\tens}{\mathbin{\otimes}}
\newcommand{\intersect}{\mathinner{\cap}}
\newcommand{\SO}{\mathrm{SO}}

\newcommand{\Spin}{\mathrm{Spin}}
\newcommand{\ie}{\textit{i.e.}}
\newcommand{\NE}{\mathcal{N\!E}}
\begin{document}
\title[Kato inequalities and conformal weights]{Refined Kato
inequalities and conformal weights in Riemannian geometry}
\author{David M. J. Calderbank} \address{D.M.J.C.: Department of
Mathematics and Statistics\\ University of Edinburgh\\ King's
Buildings, Mayfield Road\\ Edinburgh EH9 3JZ\\ Scotland.}
\email{davidmjc@maths.ed.ac.uk}
\author{Paul Gauduchon}
\address{P.G.: Centre de Math{\'e}matiques\\ Ecole Polytechnique\\
CNRS UMR 7640\\ F-91128 Palaiseau\\ France}
\email{pg@math.polytechnique.fr}
\author{Marc Herzlich}
\address{M.H.: D{\'e}partement de Math{\'e}matiques\\ Universit{\'e}
Montpellier II\\ CNRS UPRESA 5030\\ F-34095 Montpellier\\ France}
\email{herzlich@math.univ-montp2.fr}
\date{11 August 1999}
\begin{abstract}
We establish refinements of the classical Kato inequality for sections
of a vector bundle which lie in the kernel of a natural injectively
elliptic first-order linear differential operator. Our main result is
a general expression which gives the value of the constants appearing
in the refined inequalities. These constants are shown to be optimal
and are computed explicitly in most practical cases.
\end{abstract}
\maketitle

\section*{Introduction}

The Kato inequality is an elementary and well-known estimate in Riemannian
geometry, which has proved to be a powerful technique for linking 
vector-valued and scalar-valued problems in analysis on manifolds
\cite{berard-bochner,bordoni-mathann,bordoni-tsg,h-s-u,meyer-hilbert,
uhlenbeck-removable}. Its content may be stated
as follows: for any section $\xi$ of any Riemannian (or Hermitian) 
vector bundle $E$ endowed with a metric connection $\nabla$ 
over a Riemannian manifold $(M,g)$, and at any point where $\xi$
does not vanish,
\begin{equation}\label{kato1}\bigl|d|\xi|\bigr| \leq |\nabla\xi|.
\end{equation}
This estimate is easily obtained by applying the Schwarz inequality to the
right hand side of the trivial identity: $d\bigl(|\xi|^2\bigr) =
2\ip{\nabla\xi,\xi}$. Hence equality is achieved at a given point $x$ if and
only if $\nabla \xi$ is a multiple of $\xi$ at $x$, \ie, if and only if there
is a $1$-form $\alpha$ such that
\begin{equation}\label{equality} \nabla\xi = \alpha \tens \xi.
\end{equation}

The present work is motivated by circumstances in which more subtle versions
of the Kato inequality appear. Examples include: the treatment of the
Bernstein problem for minimal hypersurfaces in $\R^n$ by R.~Schoen,
L.~Simon and S.~T.~Yau \cite{ssy}, where it is shown that the second
fundamental form $h$ of any minimal immersion satisfies
\begin{equation}\label{kato-ssy} 
\bigl|d |h|\bigr| \leq \sqrt{\frac{n}{n+2}} |\nabla h|,
\end{equation} 
(see also \cite{berard-simons}); the study by S.~Bando, A.~Kasue and
H.~Nakajima of Ricci flat and asymptotically flat manifolds \cite{bkn}, where
a key role is played by the inequality
\begin{equation}\label{kato-bkn}  
\bigl|d |W|\bigr| \leq \sqrt{\frac{n-1}{n + 1}} |\nabla W|
\end{equation}
for the Weyl curvature $W$ of any Einstein metric; and the proof
given by J.~Rade of the classical decay at infinity of any Yang-Mills
field $F$ on $\R^4$ \cite{rade-twonew}, which relies on the estimate
\begin{equation}\label{kato-rade}  
\bigl|d |F|\bigr| \leq \sqrt{\frac23} |\nabla F|.
\end{equation}
Other examples may be found in the work of S.~T.~Yau on the Calabi
conjecture \cite{yau-calabiconj}, or more recently in work of P.~Feehan
\cite{feehan-decay} and of M.~Gursky and C.~LeBrun
\cite{gursky-lebrun-einstein}. For a survey of these techniques,
see also~\cite{nakajima-yautrick}.

In all of these examples, the classical Kato inequality~\eqref{kato1} is
insufficient to obtain the desired results. Moreover, the knowledge of the
best constant involved between the two terms of the inequality seems to be a
key element of all the proofs. For instance, in the case of Yang-Mills
fields on $\R^4$, the classical Kato inequality~\eqref{kato1} gives only
the decay estimate $|F| = O(r^{-2})$
at infinity, whereas the (optimal) refined inequality \eqref{kato-rade}
yields the expected $|F| = O(r^{-4})$
and thus paves the way for proving that any finite energy Yang-Mills 
field on flat space is induced from one on the sphere. 

These examples suggest that it is an interesting question to determine when
such a refined Kato inequality may occur and to compute its optimal
constant. A convincing explanation of the principle underlying this
phenomenon was first provided by J.P.~Bourguignon in \cite{jpb-magic}. He
remarked that in all the cases quoted above, the sections under
consideration are solutions of a natural linear first-order injectively
elliptic system, and that in such a situation, equality cannot occur
in~\eqref{kato1} except at points where $\nabla\xi=0$. To see this, suppose
that equality is achieved (at a point) by a solution $\xi$ of a such an
elliptic system. At that point, $\nabla \xi = \alpha\tens\xi$ for some
$1$-form $\alpha$. Now a natural first-order linear differential operator
may be written as $\Pi\circ\nabla$, where $\Pi$ is a projection onto a
(natural) subbundle of $T^*M\tens E$. Hence $\Pi(\alpha\times\xi)$
vanishes and so, by ellipticity, $\alpha\tens\xi$ vanishes.

Hence it is reasonable to expect that a refined Kato constant might appear in
this situation, \ie, that there should exist a constant $k_P<1$,
depending only on the choice of elliptic operator $P$, such that
\begin{equation}
\label{katorefined}\bigl|d |\xi|\bigr| \leq k_P |\nabla \xi|
\end{equation}
if $\xi$ lies in the kernel of $P$.

In this paper we attack the task of establishing explicitly the existence of
refined Kato constants for the injectively elliptic linear first-order
operators naturally defined on bundles associated to a Riemannian (spin)
manifold by an irreducible representation of the special orthogonal group
$\SO(n)$ or its nontrivial double-cover $\Spin(n)$. We devise a systematic
method to obtain the values of the refined constants $k_P$ and we compute
the constants explicitly in a large number of cases. We
express the constants in terms of the \emphdef{conformal weights} of
\emphdef{generalized gradients} (those operators given by projection on an
\emph{irreducible} component of the tensor product above) which are numbers
canonically attached to any such operator, and which can be easily computed
from representation theoretic data (see section 2 for details). As a
by-product of our approach, we obtain a number of representation-theoretic
formulae, relating conformal weights to higher Casimirs of $\lie{so}(n)$,
some of which appear to be new.

The structure of the paper is as follows. In the first section, we present the
basic definitions and strategy that will be followed to obtain the Kato
constants. Then, in section 2, we review the representation-theoretic
background that will be needed for our study. We do this in part for the
benefit of the reader with a limited knowledge of representation theory, but
also to set up some notation, and to demonstrate that the conformal weights
used in the sequel are easy to compute. Most importantly, we discuss the
question of which first order natural operators are injectively elliptic. This
question has been settled by Branson \cite{branson-elliptic}, whose result we
restate in the notation of this paper.

Before developing the main machinery, we use some elementary computations to
give the Kato constants when the number of irreducible components of
$T^*M\times E$ is $N=2$. Although this is entirely straightforward, the
results are sufficient to obtain a new proof of the Hijazi inequality in spin
geometry, which we sketch. For more complicated representations, we need more
tools, which we develop in section 4. Building on work of Perelomov and Popov
\cite{PP}, and also more recent ideas of Diemer and Weingart (personal
communication), we study higher Casimir elements in the universal enveloping
algebra of $\lie{so}(n)$ and obtain formulae relating them to conformal
weights. The main result in this direction is Theorem~\ref{maincas}. We use
this in sections 5 to prove our main theorem, which reduces the search for
Kato constants to linear programming. Section 6 gives some explicit constants
for $N$ odd, whereas section 7 deals with the case that $N$ is even.  In each
we give the Kato constants for a large number of operators and we detail the
precise values for $N=3$ and $N=4$. We also deal with the sharpness of our
inequalities by giving the (algebraic) equality case. Finally, as an appendix,
we present tables listing all of the Kato constants in dimensions $3$ and $4$.

\smallskip

{\sc Acknowledgements}. During the course of this work it became clear that
there is a close relationship between Kato constants and the spectral results
of Branson~\cite{branson-elliptic}. Following the presentation of an early
version of our results at a meeting in Luminy, Tom Branson has clarified this
relationship very nicely~\cite{branson-kato} and independently obtained
general minimization formula for the Kato constants. We are very grateful to
him for sharing with us his results. The formula that follows from our methods
is slightly different from his and does not cover one special case. We present
it in a similar way to permit easy comparison.

We are also deeply indebted to Tammo Diemer and Gregor Weingart for informing
us of their recent work, which plays a crucial role in our approach. Finally
we thank Christian B{\"a}r and Andrei Moroianu for the application of refined
Kato inequalities to Hijazi's inequality.

\section{Strategy}

We consider an irreducible natural vector bundle $E$ over a Riemannian (spin)
manifold $(M,g)$ of dimension $n$ with scalar product $\ip{.\,,.}$ and a
metric connection $\nabla$. By assumption, $E$ is attached to an irreducible
representation $\lambda$ of $\SO(n)$ or $\Spin(n)$ on a vector space $V$. If
$\tau$ is the standard representation on $\R^n$, then the (real) tensor
product $\tau\tens\lambda$ splits in $N$ irreducible components as
\begin{align*}
\tau\tens\lambda &= \Dsum_{j=1}^N \mu^{(j)}\\
\R^n\tens V &= \Dsum_{j=1}^N W_j.
\end{align*}
This induces a decomposition of $T^*M\tens E$ into irreducible subbundles
$F_j$ associated to the representations $\mu^{(j)}$. Projection on the
$j$th summand (of $\R^n\tens V$ or $T^*M\tens E$) will be denoted $\Pi_j$.

Following \cite{fegan,pg-pisa,hitchin-linearfields}, we can describe this
decomposition in terms of the equivariant endomorphism
$B\colon\R^n\tens V\to\R^n\tens V$ defined by
\begin{equation}\label{Boperator} 
B(\alpha \tens v)=\sum_{i=1}^n e_i\tens d\lambda(e_i\wedge\alpha) v,
\end{equation}
where $e_1,\ldots e_n$ is an orthonormal basis of $\R^n$ and $d\lambda$ is
the representation of $\lie{so}(n)$ induced by $\lambda$.
\begin{notn} For a linear map $T\colon\R^n\tens V\to\R^n\tens V$
we write $\alpha\tens\beta\mapsto T_{\alpha\tens\beta}$ for
the unique linear map $\R^n\tens\R^n\to\End(V)$ satisfying
\begin{equation}
T(\alpha\tens v)=\sum_{i=1}^n e_i\tens T_{e_i\tens\alpha}(v).
\end{equation}
Note that $(S\circ T)_{\alpha\tens\beta}
=\sum_{i=1}^nS_{\alpha\tens e_i}\circ T_{e_i\tens\beta}.$
\end{notn}
Note $B_{\alpha\tens\beta}=d\lambda(\alpha\wedge\beta)$ is a skew endomorphism
of $V$ which is skew in $\alpha\tens\beta$, and that $B$ itself is
symmetric. Therefore, the eigenvalues of $B$ are real and so, by Schur's
lemma, on the irreducible summands $W_j$, it acts by scalar multiples $w_j$ of
the identity, called \emphdef{conformal weights}. The conformal weights are
all distinct, except in the case that $V$ is an representation of $\SO(n)$
such that $\R^n\tens V$ contains two irreducible components whose sum is an
irreducible representation of $O(n)$. Therefore, apart from this exceptional
situation, the decomposition of $\R^n\tens V$ into irreducibles corresponds
precisely to its eigenspace decomposition under $B$. We shall adopt the
convention that irreducible representations of $O(n)$ in $\R^n\tens V$ will
\emph{not be split} under $\SO(n)$, so that the conformal weights $w_j$ of
$W_j$ are \emph{always distinct}. Henceforth, therefore, $W_j$ will denote the
eigenspaces of $B$ arranged so that the conformal weights $w_j$ are (strictly)
decreasing, and $N$ will denote the number of eigenspaces, \ie, the number of
(distinct) conformal weights.

The origin of this terminology is the following fact~\cite{fegan,pg-pisa}:
when the connection $\nabla$ on $E$ is induced by the Levi-Civita connection
of $(M,g)$, the natural first order operators $P_j = \Pi_j\circ\nabla$,
sometimes called \emphdef{generalized gradients}, are conformally invariant
with conformal weight $w_j$.

The operators of interest in this paper are the first order linear
differential operators $P_I:= \sum_{i\in I}\Pi_i\circ\nabla$ acting on
sections of $E$, where $I$ is a subset of $\{1,\ldots N\}$.  Such operators
are called \emphdef{Stein-Weiss operators}~\cite{stein-weiss}. The operator
$P_I$ is said to be (injectively, \ie, possibly overdetermined)
\emphdef{elliptic} iff it symbol $\Pi_I:=\sum_{i\in I} \Pi_i$ does not vanish
on any nonzero decomposable elements $\alpha\tens v$ of the tensor product
$\R^n\tens V$. Note that $P_I$ is (injectively) elliptic if and only if
$P_I^*\circ P_I$ is elliptic in the usual sense,

We could consider, more generally, the operators $\sum_{i\in I}a_i P_i$ for
any nonzero coefficients $a_i$: such an operator will be elliptic iff $P_I$
is, and the methods of this paper can be adapted to apply to this
situation. Also note that throughout the paper, $\nabla$ can be an
\emph{arbitrary} metric connection on $E$, \ie, it need not be induced by the
Levi-Civita connection of $M$.

\medbreak
We shall obtain refined Kato inequalities from refined Schwarz
inequalities of the form
\begin{equation}\label{opnorm}
\frac{|\ip{\Phi,v}|}{|v|}\leq k|\Phi|,
\end{equation}
where $\Phi\in\R^n\tens V$ and $v\in V$. For $k=1$, this holds for any $\Phi$
and nonzero $v$, with equality if $\Phi=\alpha\tens v$ for some
$\alpha\in\R^n$. Recall that the classical Kato inequality~\eqref{kato1} is
obtained from this by lifting it to the associated bundles and putting
$v=\xi$, $\Phi=\nabla\xi$ for a section $\xi$ of $E$. If $\xi$ lies in the
kernel of the operator $P_I$ then $\nabla\xi$ is a section of
$\ker\Pi_I=W_{\widehat I}$, where $\widehat I$ is the complement of $I$ in $\{1,\ldots
N\}$ and $W_{\widehat I}$ denotes the image of $\Pi_{\widehat I}$. Hence to obtain a
Kato inequality for the operator $P_I$, we only need an estimate of the
form~\eqref{opnorm} for $\Phi\in W_{\widehat I}$ and $v\in V$.  The supremum, over
all nonzero $v$, of the left hand side of~\eqref{opnorm} is the
\emphdef{operator norm} $|\Phi|_{op}$ of $\ip{\Phi,.}$, viewed as a linear map
from $V$ to $\R^n$.  Now observe that for any $\Phi\in W_{\widehat I}$, we have:
\begin{equation*}\begin{split}
|\Phi|_{op}&=\sup_{|v|=1} |\ip{\Phi,v}|
= \sup_{|\alpha|=|v|=1} |\ip{\Phi,\alpha\tens v}|
=\sup_{|\alpha|=|v|=1}|\ip{\Phi,\Pi_{\widehat I}(\alpha\tens v)}\\
&\leq\biggl(\sup_{|\alpha|=|v|=1}|\Pi_{\widehat I}(\alpha\tens v)|\biggr)|\Phi|.
\end{split}\end{equation*}
This gives a refined Schwarz inequality with
$k=\sup\limits_{|\alpha|=|v|=1}|\Pi_{\widehat I}(\alpha\tens v)|$:
\begin{equation*}\begin{split}
\frac{|\ip{\Phi,v}|}{|v|}
&=\frac{|\ip{\Phi,\alpha_0\tens v}|}{|v|}
=\frac{|\ip{\Phi,\Pi_{\widehat I}(\alpha_0\tens v)}|}{|v|}\\
&\leq \frac{|\Pi_{\widehat I}(\alpha_0\tens v)|}{|v|}|\Phi|
\leq\biggl(\sup_{|\alpha|=|v|=1}|\Pi_{\widehat I}(\alpha\tens v)|\biggr)|\Phi|,
\end{split}\end{equation*}
where $\alpha_0$ is any unit $1$-form such that $\ip{\Phi,v}=c\alpha_0$
for some $c\in\R$.

We therefore have the following Ansatz which reduces the search for refined
Kato inequalities to a purely algebraic problem.
\begin{anz}\label{ansatz}
Consider the operator $P_I$ on the natural vector bundle $E$ over $(M,g)$.
Then, for any section $\xi$ on the kernel of $P_I$, and at any point
where $\xi$ does not vanish, we have:
\begin{equation*}
\bigl|d|\xi|\bigr| \leq k_I  |\nabla\xi|,
\end{equation*}
where the constant $k_I$ is defined by
\begin{equation*}
k_I = \sup_{|\alpha|=|v|=1} |\Pi_{\widehat I}(\alpha\tens v)|.
\end{equation*}
Furthermore equality holds at a point if and only if
$\nabla\xi=\Pi_{\widehat I}(\alpha\tens\xi)$ for a $1$-form $\alpha$
at that point such that $|\Pi_{\widehat I}(\alpha\tens\xi)|=k_I|\alpha\tens\xi|$.
\end{anz}
\begin{remk} Equality holds in this Kato inequality if and only if it
holds in the refined Schwarz inequality with $v=\xi$, $\Phi=\nabla\xi$. Hence
the above Ansatz is algebraically sharp: the supremum $\sup_{|\alpha|=|v|=1}
|\Pi_{\widehat I}(\alpha\tens v)|$ is attained by compactness. We also
deduce that the Kato inequality is sharp in the flat case: equality
is attained by a suitable chosen affine solution of $P_I\xi=0$.
\end{remk}

\noindent In order to turn this Ansatz into a useful result, we must:
\begin{enumerate}
\item Find when $P_I$ is elliptic.
\item Show that when $P_I$ is elliptic, $k_I$ is less than one.
\item Give a formula for $k_I$ in terms of easily computable data.
\item Obtain a more explicit description of the equality case.
\end{enumerate}

The first question has been answered by T. Branson~\cite{branson-elliptic}. We
shall discuss his result at the end of the next section. Also in that section
we shall give a more explicit description of the operators and representations
involved, together with the associated conformal weights. The conformal
weights are easy to compute and so our guiding philosophy will be: \emph{find
$k_I$ in terms of the conformal weights}.

Since $k_I = \sup_{|\alpha|=|v|=1} |\Pi_{\widehat I}(\alpha\tens v)|$ and
\begin{equation*}
|\Pi_{\widehat I}(\alpha\tens v)|^2
=\sum_{j\in\widehat I}|\Pi_j(\alpha\tens v)|^2,
\end{equation*}
a key step in our task is to find a convenient formula for
$|\Pi_j(\alpha\tens v)|$ for each $j=1,\ldots N$. 

To do this, note that $\Pi_j$ is the projection onto an eigenspace of $B$,
and so Lagrange interpolation gives the standard formulae:
\begin{equation}
\Pi_j=\prod_{k\neq j} \frac{B-w_k\iden}{w_j-w_k}
=\frac{\displaystyle\sum_{k=0}^{N-1}
w_j^{N-1-k}\sum_{\ell=0}^k(-1)^\ell\sigma_\ell(w)B^{k-\ell}}
{\displaystyle\prod_{k\neq j}(w_j-w_k)},
\end{equation}
where $\sigma_i(w)$ denotes the $i$th elementary symmetric function in
the eigenvalues $w_j$. We define $A_k$ to be the operators
\begin{equation}
A_k=\sum_{\ell=0}^k(-1)^\ell\sigma_\ell(w)B^{k-\ell}
\end{equation}
appearing in this formula, which are manifestly symmetric in the conformal
weights. Using these operators, we have:
\begin{equation}\label{proj}
|\Pi_j(\alpha\tens v)|^2=\ip{\Pi_j(\alpha\tens v),\alpha\tens v}
=\frac{\displaystyle\sum_{k=0}^{N-1} w_j^{N-1-k}
\ip{A_k(\alpha\tens v),\alpha\tens v}}{\displaystyle\prod_{k\neq j}(w_j-w_k)}.
\end{equation}

This formula for the $N$ quantities $|\Pi_j(\alpha\tens v)|$ in terms of the
$N$ quantities $q_k=\ip{A_k(\alpha\tens v),\alpha\tens v}$ lies at the heart
of our method. Note first that $A_0=1$, and so $q_0=|\alpha\tens v|^2$, which
we set equal to $1$. Secondly, the formula~\eqref{Boperator} for $B$ implies
that
\begin{equation}\label{Bzero}
\ip{B(\alpha\tens v),\alpha\tens v}=0,\qquad\forall\alpha \in \R^n, v\in V.
\end{equation}
Hence $q_1$ is also computable. These two observations alone will allow us to
find the Kato constants for $N\leq4$. For larger $N$ we shall need to obtain
more information about the operators $A_k$.

We shall find that approximately half of the $q_k$'s can be eliminated. The
remainder can then be estimated from above and below using the non-negativity
of $|\Pi_j(\alpha\tens v)|$. These bounds can in turn be used to estimate
$|\Pi_{\widehat I}(\alpha\tens v)|$.

\section{Representation theoretic background}\label{repback}

The description of representations of the special orthogonal group $\SO(n)$,
or its Lie algebra $\lie{so}(n)$ differs slightly according to the parity of
$n$. We write $n=2m$ if $n$ is even and $n=2m+1$ if $n$ is odd; $m$ is then
the \emphdef{rank} of $\lie{so}(n)$.

We fix an oriented orthonormal basis $(e_1, \ldots e_n)$ of $\R^n$, so that
$e_i\wedge e_j$ (for $i<j$) is a basis of the Lie algebra $\lie{so}(n)$,
identified with $\Lambda^2\R^n$. We also fix a Cartan subalgebra $\lie{h}$ of
$\lie{so}(n)$ by the basis $E_1 = e_1 \wedge e_2, \ldots E_m = e_{2m - 1}
\wedge e_{2m}$, and denote the dual basis of $\lie{h}^*$ by $(\eps_1, \ldots
\eps_m)$. We normalize the Killing form so that this basis is
orthonormal. For further information on this, and the following, see
\cite{fulton-harris,salamon-holonomy,samelson}.

An irreducible representation of $\lie{so}(n)$ will be identified with its
dominant weight $\lambda\in\lie{h}^*$. Roots and weights can be given by
their coordinates with respect to the orthonormal basis $\eps_i$.  Then the
weight $\lambda=(\lambda_1, \lambda_2,\ldots \lambda_m)$, whose coordinates
are all integers or all half-integers, is dominant iff
\begin{align*}
\lambda_1 \geq\cdots\geq \lambda_{m-1}\geq|\lambda_m|, \qquad&n=2m,\\
\lambda_1 \geq\cdots\geq \lambda_{m-1}\geq\lambda_m\geq0, \qquad&n=2m+1.
\end{align*}
In this notation, the standard representation $\tau$ is given by the weight
$(1,0,...0)$, the weight $\lambda=(1,1,...1,0,...0)$ (with $k$ ones)
corresponds to the $k$-form representation $\Lambda^k\R^n$, the weights
$\lambda=(1,1,...1,\pm1)$ (for $n=2m$) correspond to the selfdual and
antiselfdual $m$-forms, and the weights
$\lambda=\bigl(\tfrac12,...\tfrac12,(\pm)\tfrac12\bigr)$ correspond to the
spin or half-spin representations $\Delta_{(\pm)}$. The \emphdef{Cartan
product} of two representations is the subrepresentation $\lambda\odot\mu$ of
highest weight $\lambda+\mu$ in $\lambda\tens\mu$. If $\lambda$ and $\mu$ are
integral then $\lambda\odot\mu$ is the subrepresentation of
``alternating-free, trace-free'' tensors in $\lambda\tens\mu$; for instance,
the $k$-fold Cartan product $\odot^k\R^n$ is the representation $S^k_0\R^n$ of
totally symmetric traceless tensors, with weight $(k,0,...0)$.

Notice that we take the real form of the representations wherever possible:
in particular, when discussing elements of the tensor product
$\tau\tens\lambda$, only real elements of the standard representation will be
used, even if $\lambda$ is complex.

The decomposition of the tensor product $\tau\tens\lambda$ into
irreducibles is described by the following rule: an irreducible
representation of weight $\mu$ appears in the decomposition if and only if
\begin{enumerate}
\item $\mu = \lambda\pm\eps_j$ for some $j$, or $n=2m+1$, $\lambda_m>0$
and $\mu=\lambda$
\item $\mu$ is a dominant weight.
\end{enumerate}
Weights $\mu$ satisfying (i) will be called \emphdef{virtual} weights
associated to $\lambda$. We shall say $\mu$ is \emphdef{effective} if it also
satisfies (ii). It will be convenient to have a notation for the virtual
weights which is compatible with the outer automorphism equivalence of
representations of $\lie{so}(2m)$. We define $\mu^0=\lambda$ and
$\mu^{i,\pm}=\lambda\pm\eps_i$, unless $n=2m$, $j=m$ and $\lambda_m\neq0$, in
which case we define $\mu^{m,\pm}$ to be the virtual weights such that
$|\mu^{m,+}_m|=|\lambda_m|+1$ and $|\mu^{m,-}_m|=|\lambda_m|-1$. This notation
allows us to assume, without loss of generality, that $\lambda_m=|\lambda_m|$,
and we shall omit the modulus signs in the following.

The \emphdef{Casimir number} of a representation $\lambda$ is given by 
\begin{equation}\label{cas}
c(\lambda) =\ip{\lambda+\delta,\lambda+\delta}-\ip{\delta,\delta}
=\ip{\lambda,\lambda} + 2\ip{\lambda,\delta},
\end{equation}
where $\delta$ is the half-sum of positive roots, \ie,
$\delta_i=(n-2i)/2$.

The \emphdef{conformal weight} associated to a component $\mu$ of
$\tau\tens\lambda$ may be computed explicitly by the formula
\begin{equation}\label{cwts}
w(\mu,\lambda) = \tfrac12 (c(\mu)-c(\lambda)-c(\tau)),
\end{equation}
which continues to make sense for virtual weights. We let $w^0$ and
$w^{i,\pm}$ denote the (virtual) conformal weights of $\mu^0$ and
$\mu^{i,\pm}$. Expanding the definition of the Casimir, and applying some
Euclidean geometry in $\lie{h}^*$, we obtain the explicit formulae (assuming
$\lambda_m=|\lambda_m|$):
\begin{align}\label{explicit}
w^0&=(1-n)/2\\
w^{i,+}&=1+\lambda_i-i\\
w^{i,-}&=1-n-(\lambda_i-i).
\end{align}
These formulae show that conformal weights are simple to compute in practice,
which is one of our motivations for using them.  We note that the virtual
conformal weights $w^{i,\pm}$ satisfy
\begin{equation}
w^{1,+} > w^{2,+} > \cdots > w^{m,+} \geq w^{m,-} > \cdots > w^{2,-} > w^{1,-}
\end{equation}
with equality in the middle if and only if $n=2m$ and $\lambda_m=0$.
If $n=2m+1$ and $\lambda_m>0$, then the conformal weight $w^0$ lies
strictly between $w^{m,+}$ and $w^{m,-}$. This verifies our earlier
claim that the conformal weights are almost always distinct.

For effective weights, we remind the reader of our convention not to split
subrepresentations with the same conformal weight. This means that we write
$\tau\tens\lambda=\oplus_{j=1}^N \mu^{(j)}$, where the representations
$\mu^{(j)}$ are all irreducible, unless $n=2m$ and $\lambda_m=0$, in which one
of the components is taken to be $\mu^{m,+}\dsum\mu^{m,-}$.

In order to say which of the weights are effective (and hence, which
representations occur in $\tau\tens\lambda$), it is useful to make explicit
repetitions among the coordinates $\lambda_j$ by writing $\lambda$ in the form:
\begin{equation*}
\lambda=(k_1,\ldots k_1, k_2,\ldots k_2,\, \ldots\ldots, k_\nu,\ldots k_\nu),
\end{equation*}
with $k_1 > k_2 > \cdots > k_{\nu-1} > k_{\nu} \geq 0$. If $k_\nu\neq0$
and $n=2m$, we write
\begin{equation*}
\lambda=(k_1,\ldots k_1, k_2,\ldots k_2,\, \ldots\ldots, k_\nu,\ldots\pm k_\nu)
\end{equation*}
for the two possible signs of the last entry.  Here $\nu$ is
the number of groups of equal entries and we let $p_1$ denote the number of
$k_1$'s, $p_2-p_1$ the number of $k_2$'s, etc., so that $p_j$ is the number of
entries greater than or equal to $k_j$.

We first note that the following $2\nu-1$ weights, at least, are effective for
any representation $\lambda$, and are associated with the conformal weights
listed.
\begin{align*}
&\mu^{1,+}          &w^{1,+}&=k_1\\
&\mu^{p_1+1,+}      &w^{p_1+1,+}&=k_2-p_1\\
&\;\vdots           &&\;\;\vdots\\
&\mu^{p_{\nu-1}+1,+}&w^{p_{\nu-1}+1,+}&=k_\nu-p_{\nu-1}\\
&&&\\
&\mu^{p_{\nu-1},-}  &w^{p_{\nu-1},-}&=p_{\nu-1}-k_{\nu-1}+1-n\\
&\;\vdots           &&\;\;\vdots\\
&\mu^{p_1,-}        &w^{p_1,-}&=p_1-k_1+1-n.
\end{align*}
If $k_\nu=0$ there are no further effective weights unless $n=2m$ and
$p_{\nu-1}=m-1$, in which case $\mu^{m,\pm}$ are both effective with the
same conformal weight. Hence, by convention, if $k_\nu=0$ then $N=2\nu-1$.

If $k_\nu>0$ and $n=2m$ then $\mu^{m,-}$ is effective and $N=2\nu$. If
$k_\nu>0$ and $n=2m+1$ then $\mu^0$ is a possible target; furthermore
$\mu^{m,-}$ is effective for $k_\nu>1/2$.

We therefore see that the number of components $N$ in the decomposition
$\R^n\tens V=\oplus_{j=1}^N W_j$ is either $2\nu-1$, $2\nu$ or $2\nu+1$.

The case $N=2\nu-1$ arises when $\lambda_m=0$. The representations occuring,
in order of decreasing conformal weight are as follows.
\begin{align*}
\mu^{(1)}&=\mu^{1,+},&\mu^{(2)}&=\mu^{p_1+1,+},&\ldots
&&\mu^{(\nu-1)}&=\mu^{p_{\nu-2}+1,+},\\
\mu^{(\nu)}&=\mu^{p_{\nu-1}+1,+}\quad\textrm{or}\quad\mu^{m,+}\dsum\mu^{m,-}
,&&&&&&\\
\mu^{(\nu+1)}&=\mu^{p_{\nu-1},-},&\mu^{(\nu+2)}&=\mu^{p_{\nu-2},-},
&\ldots&&\mu^{(2\nu-1)}&=\mu^{p_1,-}.
\end{align*}

The case $N=2\nu$ arises when $n=2m+1$ and $\lambda_m=1/2$ or when
$n=2m$ and $\lambda_m\neq0$. The representations occuring, in
order of decreasing conformal weight are as follows.
\begin{align*}
\mu^{(1)}&=\mu^{1,+},&\mu^{(2)}&=\mu^{p_1+1,+},&\ldots
&&\mu^{(\nu)}&=\mu^{p_{\nu-1}+1,+},\\
\mu^{(\nu+1)}&=\mu^{m,-}\quad\textrm{or}\quad\mu^0,&
\mu^{(\nu+2)}&=\mu^{p_{\nu-1},-},&\ldots&&\mu^{(2\nu)}&=\mu^{p_1,-}.
\end{align*}

The case $N=2\nu+1$ arises when $n=2m+1$ and $\lambda_m>1/2$. The
representations occuring, in order of decreasing conformal weight are as
follows.
\begin{align*}
\mu^{(1)}&=\mu^{1,+},&\mu^{(2)}&=\mu^{p_1+1,+},&\ldots
&&\mu^{(\nu)}&=\mu^{p_{\nu-1}+1,+},\\
\mu^{(\nu+1)}&=\mu^0,&&&&&&\\
\mu^{(\nu+2)}&=\mu^{m,-},&\mu^{(\nu+3)}&=\mu^{p_{\nu-1},-},
&\ldots&&\mu^{(2\nu+1)}&=\mu^{p_1,-}.
\end{align*}
Note that for ``most'' representations (e.g., if $\lambda_m\neq0$) $N$ and $n$
have the same parity. Indeed, if $\lambda_1>\lambda_2>\cdots >|\lambda_m|>0$
we see that $N=n$. However, the representations arising in practice are not at
all generic: $N$ is usually very small.

We are now in a position to describe T. Branson's classification of the
elliptic operators~\cite{branson-elliptic}.  Firstly, note that if $J$ is a
subset of $I$ such that $P_J$ is elliptic, then $P_I$ is elliptic. Hence it
suffices to find the \emphdef{minimal elliptic operators} $P_I$, \ie, the
elliptic $P_I$ such that $P_J$ is not elliptic for any proper subset $J$ of
$I$.

\begin{thm}[Branson \cite{branson-elliptic}]\label{ellipresult} Let
$\lambda$ be an irreducible representation of $\SO(n)$ or $\Spin(n)$. Then the
minimal elliptic operators associated to $\lambda$ are either elementary or
the sum of two elementary operators.  The elementary elliptic operators are:
\begin{enumerate}
\item $P_1$ with target $\mu^{1,+}$.
\item For $N=2\nu:$\quad $P_{\nu+1}$ with target $\mu^{m,-}$ or $\mu^0$.
\item For $N=2\nu+1:$\quad $P_{\nu+1}$ with target $\mu^0$,
provided $\lambda$ is properly half-integral.
\end{enumerate}
The other minimal elliptic operators are:
\begin{enumerate}\setcounter{enumi}3
\item $P_{\{j,N+2-j\}}$ with target $\mu^{p_{j-1}+1,+}\dsum\mu^{p_{j-1},-}$ or
$\mu^{m,+}\dsum\mu^{m,-}\dsum\mu^{m-1,-}$ for all $j\in\{2,...\nu\}$.
\textup(For $N=2\nu-1$, $j=\nu$ and $p_{\nu-1}=m-1$, $P_{\nu,\nu+1}$ is
obtained by combining the operators with targets
$\mu^{m,\pm}\dsum\mu^{m-1,-}$, which are both elliptic.\textup)
\item For $N=2\nu+1:$\quad $P_{\{\nu+1,\nu+2\}}$ with target
$\mu^0\dsum\mu^{m,-}$, provided $\lambda$ is integral.
\end{enumerate}
\end{thm}
Notice that the subsets of $N$ corresponding to the minimal elliptic operators
partition $N$ (where we combine the operators with targets
$\mu^{m,\pm}\dsum\mu^{m-1,-}$), unless $N=2\nu+1$ and $\lambda$ is properly
half-integral, in which case there is one ``useless'' operator $P_{\nu+2}$.
This means that there are non-elliptic operators with relatively large
targets. Indeed, the above theorem may equivalently be viewed as a description
of the maximal non-elliptic operators. These play an important role in our
later work, so we shall describe them explicitly here.

\begin{defn} Let $\NE$ denote the set of subsets of $\{1,\ldots N\}$
whose elements are obtained by choosing exactly one index in each of the sets
$\{j,N+2-j\}$ for each $j$ with $2\leq j\leq\nu$ if $N=2\nu-1, 2\nu$ (giving
$2^{\nu-1}$ elements) and for each $j$ with $2\leq j\leq\nu+1$ if $N=2\nu+1$
(giving $2^\nu$ elements).
\end{defn}
Branson's theorem implies that the set $\NE$ is precisely the set of subsets
of $\{1,\ldots N\}$ corresponding to the maximal non-elliptic operators,
unless $N=2\nu+1$ and $\lambda$ is properly half-integral, in which case the
maximal non-elliptic operators correspond to the elements of $\NE$ which do
not contain $\nu+1$. This last case will cause us problems because there are
not enough non-elliptic subsets.

Branson proves Theorem~\ref{ellipresult} by reducing the problem to the study
of the spectrum of the operator on the sphere $M=S^n$, which he computes by
applying powerful techniques from harmonic analysis. For the benefit of the
reader not familiar with these global techniques, we remark that there are
some cases in which ellipticity or non-ellipticity can be established by
elementary local arguments.

Since ellipticity depends only on the symbol $\Pi_I$ on $\R^n\tens V$ and
since $\SO(n)$ is transitive on the unit sphere in $\R^n$, it follows that
$P_I$ is elliptic if and only if the linear map $v\to\Pi_I(e_n\tens v)$ is
injective (for a fixed unit vector $e_n$).

First note that this map is $\SO(n-1)$-equivariant and so we have the
following necessary (but not sufficient) condition for ellipticity.

\begin{lemma}\label{nonelliptic} $P_I$ cannot be elliptic unless every
subrepresentation of $V$ under the group $\SO(n-1)$ occurs as a
subrepresentation of $W_j$ for some $j\in I$.
\end{lemma}

To use this lemma, one must apply the standard branching rule branching rule
for restricting a representation of $\SO(n)$ to $\SO(n-1)$---see, for
example~\cite[page 426]{fulton-harris}. For $N=2\nu-1$ and $N=2\nu$ it is
straightforward to verify the non-ellipticity of the maximal non-elliptic
operators and hence obtain most of the non-ellipticity results in Branson's
theorem. For $N=2\nu+1$ this naive method does not cover all the cases:
$P_{\nu+1}$ is \emph{not} elliptic if $\lambda$ is an integral weight, even
though $\lambda$ itself is the target representation.

Secondly, note the following sufficient (but not necessary) condition
for ellipticity.

\begin{lemma} If the space of local solutions of $P_I$ on $\R^n$ is finite
dimensional, then $P_I$ is elliptic.
\proofof{lemma} If $P_I$ is not elliptic then for some $v\in V$,
$e_n\tens v$ belongs to $\ker\Pi_I\leq \R^n\tens V$.
Now consider the operator $P_I$ on $\R^n$ (with respect to the
trivial connection on $E$). If $L_v$ denotes the line subbundle of $E$
corresponding to the span of $v\in V$ then any section of $L_v$ which
is independent of $x_1,...x_{n-1}$ belongs the kernel of $P_I$.
Hence the kernel of $P_I$ is infinite dimensional on $\R^n$.
\end{proof}

As observed (for instance) in~\cite{KOPWZ}, this second lemma shows that the
\emphdef{highest gradient} is always elliptic. This is the operator $P_1$
with the highest conformal weight $w_1$ whose target $\mu^{(1)}$ is the highest
weight subrepresentation of $\tau\tens\lambda$. We shall also refer to $P_1$
as the \emphdef{Penrose} or \emphdef{twistor operator}, since it reduces to
the usual Penrose twistor operator if one views the representation $\lambda$
as a subrepresentation of a tensor product of spinor representations. The
kernel of a twistor operator on $S^n$ (or any simply connected open subset)
is well-known to be a finite dimensional representation space for
$\SO(n+1,1)$: the twistor operator is the first operator in the
Bernstein-Gelfand-Gelfand resolution of this representation (see for
instance~\cite{be}).

Finally in this section, we recall the following ellipticity result:
\begin{prop}\tcite{pg-pisa} $P_I$ is elliptic in either of the following
cases:
\begin{enumerate}
\item $I$ contains all $j$ with $w_j\geq0$
\item $I$ contains all $j$ with $w_j\leq0$.
\end{enumerate}
\end{prop}
These operators are of special interest because there
is a simple Weitzenb\"ock formula relating them~\cite{pg-pisa}.

\section{Refined Kato inequalities with $N=2$}

The case $N=2$ often arises in spin geometry and in two and four dimensional
differential geometry. It occurs in the following two cases:
\begin{enumerate}
\item When the dimension $n$ is even, $\lambda = (k,\ldots k, \pm k)$ with
$k$ an arbitrary integer or half-integer, \ie, $V = \odot^{2k} \Delta_+$ or $V
= \odot^{2k} \Delta_-$. Therefore the bundle $E$ is either $\odot
^k\Lambda^m_{\pm}M$ or, if $M$ is spin, $\odot^{k-\frac12}\Lambda^m
_{\pm}M \odot \Sigma^{\pm}$ ($\Sigma^{\pm}$ denote positive and negative
spinor bundles of $M$); one thus gets $w_1 = k > w_2 = 1 - \frac{n}2  - k$.
\item When the dimension $n$ is odd, $\lambda = (\frac12,
\ldots  \frac12)$, \ie, $V = \Delta$, $E$ is the spinor bundle
$\Sigma$ and $w_1 = \frac12 \ > \  w_2 =  \frac{1 - n}2$. 
\end{enumerate}
\noindent Note that the operators $P_1$ and $P_2$ are both elliptic.
\begin{thm}\label{thkatoN2} Let $E$ be associated to a representation
$\lambda$ with $N=2$.
\begin{enumerate}
\item For any nonvanishing section $\xi$ of $E$ in the kernel of the twistor
operator $P_1$,
\begin{equation} \label{kato2i} 
\bigl|d|\xi|\bigr| \leq \sqrt{\frac{w_1}{w_1 - w_2}}  |\nabla \xi| 
= \sqrt{\frac{k}{2k + \frac n2 - 1}}  |\nabla \xi| 
\end{equation}
with equality if and only if, for some $1$-form $\alpha$,
\begin{equation*} 
  \nabla \xi = \Pi _2 (\alpha\tens\xi).
\end{equation*}
\item For any section $\xi$ of $E$ in the kernel of $P_2$,     
\begin{equation} \label{kato2ii}
\bigl|d|\xi|\bigr|  \leq \sqrt{\frac{-w_2}{w_1 - w_2}} |\nabla \xi|
= \sqrt{\frac{k + \frac n2 - 1}{2k +\frac n2 - 1}}  |\nabla \xi| 
\end{equation}
with equality if and only if, for some $1$-form $\alpha$,
\begin{equation*} 
  \nabla \xi = \Pi _1 (\alpha \tens \xi). 
\end{equation*}
\end{enumerate}
\proofof{thm} From the Ansatz~\ref{ansatz}, we have to estimate the norms of
$\Pi_j(\alpha\tens v)$ for $j=1,2$. The crucial ingredient here is
equation~\eqref{Bzero}, which gives the following system of equations
for the components of a unit length vector $\alpha\tens v$ in $\R^n\tens V$:
\begin{equation}\label{syst2}\begin{split}
& |\Pi _1 (\alpha \tens v)|^2 + |\Pi_2 (\alpha \tens v)|^2 = 1,\\
& w_1 |\Pi _1 (\alpha \tens v)|^2 + w_2 |\Pi_2 (\alpha \tens v)|^2 = 0.
\end{split}\end{equation}
The solution is a special case of equation~\eqref{proj}:
\begin{equation} \label{decpi}
|\Pi_1 (\alpha \tens v)|^2=\frac{w_2}{w_2 - w_1},
\hspace{0.5cm} |\Pi _2 (\alpha \tens v)|^2 = \frac{w_1}{w_1 - w_2}
\end{equation}
and moreover this is valid for {\it any} choice of unit $\alpha$ and $v$.
These formulae easily yield the refined Kato inequalities and
their equality cases.
\end{proof}
\begin{remk}\label{weitz} The calculations above also yield some 
(possibly not optimal) refined Kato inequalities for $N$ arbitrary and
operators
\begin{equation*}
P_+ = \sum_{w_j >0} P_j \ \textrm{ or } \   P_- = \sum_{w_j<0} P_j
\end{equation*}
(for simplicity's sake, we consider here only the case when conformal weights
do not vanish). The reasoning for $P_+$ relies on the system of
equations
\begin{equation}\begin{cases} 
|\Pi _+ (\alpha \tens v)|^2 + |\Pi_- (\alpha \tens v)|^2 = 1,\\ 
 w_1  \, |\Pi _+ (\alpha \tens v)|^2 
  + w^{max}_{< 0}  \, |\Pi_- (\alpha \tens v)|^2 \geq  0
\end{cases}\ \ \textrm{ with }\  w^{max}_{< 0}  = \max_{w_j<0} w_j 
\end{equation}
with $\Pi_{\pm}$ the projections associated to both operators.
One easily gets the refined Kato inequality 
\begin{equation} 
\bigl|d|\xi|\bigr| \leq\sqrt{\frac{w_1}{w_1 - w^{max}_{< 0}}}\,|\nabla \xi| 
\end{equation}
for any section $\xi$ in the kernel of $P_+$ and similarly 
\begin{equation} 
\bigl|d|\xi|\bigr| \leq \sqrt{\frac{w_N}{w_N - w^{min}_{>0}}}\,|\nabla \xi|,  
\quad\textrm{with}\quad w^{min}_{>0} = \min_{w_j>0} w_j
\end{equation}
for any section $\xi$ in the kernel of $P_-$.
\end{remk}

\begin{remk}\label{katohijazi}  As an application of these results,
we give a new proof of the {\it Hijazi inequality} in spin geometry
relating the first eigenvalue of the Dirac operator on a Riemannian spin
manifold to the first eigenvalue of its conformal Laplacian. This application
is due to Christian B{\"a}r and Andrei Moroianu 
(private communication), and we thank them for their
permission to reproduce it in this work.

\begin{prop}[Hijazi \cite{hijazi-conformal}] 
Let $(M,g)$ be a compact Riemannian spin manifold of dimension 
$n\geq 3$. Then the first eigenvalue $\lambda_1$ of the Dirac operator and 
the first eigenvalue $\mu_1$ of the conformal Laplacian
$4\frac{n-1}{n-2}\Delta + \Scal$ satisfy:
\begin{equation}\label{hij}
\lambda_1^2 \geq \frac{n}{4(n-1)}\, \mu_1.
\end{equation}
\proofof{prop} If $\psi$ is an eigenspinor with eigenvalue $\lambda$, then
$\psi$ lies in the kernel of the Dirac operator given by the Friedrich
connection $\tilde\nabla_X\psi = \nabla_X\psi + (\lambda/n)
X\cdot\psi$, which is a metric connection on spinors. Hence we have the
following refined Kato inequality for $\psi$, wherever it is nonzero:
\begin{equation}\label{kad}
\bigl|d|\psi|\bigr|^2\leq \frac{n-1}n|\tilde\nabla\psi|^2.
\end{equation}
We next consider the conformal Laplacian of $|\psi|^{2\alpha}$ where
$\alpha=\frac{n-2}{2(n-1)}$: the conformal Laplacian is invariant on scalars
of weight $\frac{2-n}2$ and so this power is natural in view of the conformal
weight $\frac{1-n}2$ for the Dirac operator. Using the Lichnerowicz formula
and the elementary identity $d^*d(f^\alpha)=\alpha
f^{\alpha-1}d^*df-\alpha(\alpha-1) f^{\alpha-2}|df|^2$ with $f=|\psi|^2$, we
obtain the following equalities on the open set where $\psi$ is nonzero:
\begin{align*}
\tfrac1{2\alpha}d^*d\bigl(&|\psi|^{2\alpha}\bigr)
+\tfrac14\Scal\,|\psi|^{2\alpha}-\tfrac{n-1}n\lambda^2|\psi|^{2\alpha}\\
&=\tfrac12(1-\alpha)|\psi|^{2\alpha-4}\bigl|d\bigl(|\psi|^2\bigr)\bigr|^2
+\tfrac12|\psi|^{2\alpha-2}d^*d\bigl(|\psi|^2\bigr)
+\bigl(\tfrac14\Scal-\tfrac{n-1}n\lambda^2\bigr)|\psi|^{2\alpha}\\
&=|\psi|^{2\alpha-2}\Bigl(2(1-\alpha)\bigl|d|\psi|\bigr|^2
+\ip{\nabla^*\nabla\psi,\psi}-|\nabla\psi|^2
+\tfrac14\Scal\,|\psi|^2-\tfrac{n-1}n\lambda^2|\psi|^2\Bigr)\\
&=|\psi|^{2\alpha-2}\Bigl(2(1-\alpha)\bigl|d|\psi|\bigr|^2
+\bigl(1-\tfrac{n-1}n\bigr)\lambda^2|\psi|^2-|\nabla\psi|^2\Bigr)\\
&=|\psi|^{2\alpha-2}\Bigl(\tfrac{n}{n-1}\bigl|d|\psi|\bigr|^2
-|\tilde\nabla\psi|^2\Bigr)
\end{align*}
since $|\tilde{\nabla}\psi|^2=|\nabla\psi|^2+\frac1n\lambda^2|\psi|^2$.
This is nonpositive by~\eqref{kad}. Notice that this gives a local version
of the Hijazi inequality, with equality iff $\tilde\nabla\psi$ is the
projection of $\alpha\tens\psi$ onto the kernel of Clifford
multiplication, for some $1$-form $\alpha$. If the eigenvalue $\lambda$
is nonzero, then differentiating and commuting derivatives shows in fact
that $\tilde\nabla\psi=0$. The case $\lambda=0$ is distinguished by
conformal invariance and the fundamental solutions $\psi(x)=c(x)\phi/|x|^n$
give examples with $\nabla\psi\neq0$.

In order to globalize, we consider the Rayleigh quotient for the first
eigenvalue $\mu_1$ of the conformal Laplacian:
\begin{equation*}
\mu_1 \leq \frac{\int_M 4\frac{n-1}{n-2}|d\varphi|^2+\Scal\,\varphi^2}
{\int_M \varphi^2}.
\end{equation*}
We can estimate the integral in the numerator by setting
$\varphi=|\psi|^{2\alpha}$ on the open set where $\psi$ is nonzero and
writing
\begin{align*}
4\frac{n-1}{n-2}\bigl|\bigl(d|\psi|^{2\alpha}\bigr)\bigr|^2
&+\Scal\,|\psi|^{4\alpha}\\
&=4|\psi|^{2\alpha}\Bigl(\frac{n-1}{n-2}d^*d\bigl(|\psi|^{2\alpha}\bigr)
+\frac14\Scal\,|\psi|^{2\alpha}\Bigr)
-\frac{2(n-1)}{n-2}d^*d\bigl(|\psi|^{4\alpha}\bigr)\\
&\leq\frac{4(n-1)}n\lambda^2|\psi|^{4\alpha}
-\frac{2(n-1)}{n-2}d^*d\bigl(|\psi|^{4\alpha}\bigr).
\end{align*}
Taking $\lambda=\lambda_1$, integrating over $\{x\in M:|\psi|(x)\geq\eps\}$
and letting $\eps\to0$ gives~\eqref{hij}. The equality case is also
easy to establish.
\end{proof}
\end{remk}
A similar argument can be used to provide an alternative proof the
$N=2$ vanishing theorems of Branson-Hijazi~\cite{branson-hijazi1}.

\section{Casimir numbers and conformal weights}

One way to understand the powers $B^\ell\colon\R^n\tens V\to\R^n\tens V$ of
the operator $B$ is to relate them to invariants of $V$. Let $\ptr
B^\ell=\sum_i(B^\ell)_{e_i\tens e_i}\colon V\to V$ be the \emphdef{partial
trace} of $B$ obtained by contracting over $\R^n$. Since $V$ is irreducible
and $B$ is symmetric and equivariant, this partial trace must be a scalar
multiple of the identity. The explicit expression~\eqref{Boperator} for $B$
yields the following formula:
\begin{equation}
\ptr B^{\ell}=\sum_{i_1,\ldots i_\ell}
d\lambda(e_{i_1}\wedge e_{i_2})\circ d\lambda(e_{i_2}\wedge e_{i_3})\circ
\cdots\circ d\lambda(e_{i_{\ell-1}}\wedge e_{i_\ell})\circ
d\lambda(e_{i_\ell}\wedge e_{i_1}).
\end{equation}
This is the action on $V$ of an element of the centre of the universal
enveloping algebra $\mathcal{U}(\lie{so}(n))$ called a \emphdef{higher
Casimir}, since it reduces to the Casimir element when $\ell=2$ (and
vanishes when $\ell=1$). The (scalar) action of the Casimir element on $V$ is
the Casimir number $c(\lambda)$ of $V$, and it is of some interest to compute
the \emphdef{higher Casimir numbers}.  This computation was carried out by A.
Perelomov and V. Popov in \cite{PP}, where a generating series for the higher
Casimir numbers in terms of polynomials in $\lambda$ is given.

Our aim in this section is to obtain instead relations between higher
Casimirs and conformal weights. These relations will enable us to find
a more convenient basis for the higher Casimirs in terms of certain
linear combinations of the $B^\ell$.

In fact it is more natural to work with the \emphdef{translated} operator
$\widetilde B=B+\frac{n-1}2\iden$ and its eigenvalues, the \emphdef{translated
conformal weights} $\widetilde w_j=w_j+\frac{n-1}2 =
\frac12\bigl(c(\mu^{(j)})-c(\lambda)\bigr)$. The translated virtual conformal
weights are then $\widetilde w^{i,\pm}=\frac12 \pm (\lambda_i+\frac n2-i)
=\frac12\pm x_i$ where $x=\lambda+\delta$. These translated conformal weights
are more convenient because if $\lambda_i=\lambda_{i+1}$ then
\begin{equation}\label{cancel}
\widetilde w^{i+1,+} + \widetilde w^{i,-} = 0.
\end{equation}
which is a useful cancellation property for non-effective weights.
In particular, there is the following immediate consequence, which
already suggests that (translated) conformal weights are a convenient
tool for handling Casimir numbers.

\begin{prop}\label{propid} Let $P_{\ell}$ be the polynomial on \textup(the
dual of\textup) the Cartan subalgebra defined by
\begin{equation*}
P_{\ell} (\lambda) = \sum _{i = 1}^m \Bigl(\frac12 + x_i\Bigr)^{\ell}
 + \sum_{i = 1}^m \Bigl(\frac12 - x_i\Bigr)^{\ell}
\qquad\textrm{for}\quad \ell \in \N,
\end{equation*}
where $x=\lambda+\delta$. Then:
\begin{enumerate}
\item if $N$ is odd,
\begin{equation} \label{idodd}   \sum _{j = 1}^N
  \widetilde w_j^{2k + 1}  -  \Bigl(\frac{n - 1}2\Bigr)^{2k + 1}=
P_{2k + 1}(x) - P_{2k + 1} (\delta) \qquad\forall\,  k \in \N;
\end{equation}
\item if $N$ is even
\begin{equation} \label{ideven}  \sum_{j = 1}^N
  \widetilde w_j^{2k + 1} -  \Bigl(\frac{n - 1}2\Bigr)^{2k + 1}
- \Bigl(\frac12\Bigr)^{2k + 1}
= P_{2k + 1}(x) - P_{2k + 1} (\delta) \qquad\forall\, k \in \N.
\end{equation}
\end{enumerate}
\proofof{prop} The starting point is the trivial formula
\begin{equation*}
P_{2k+1} (x)=\sum (\widetilde w^{i,\pm})^{2k+1} \quad\forall\,k\in \N,
\end{equation*}
where the summation is over all virtual weights
(it does not matter whether we include $\mu^0$ as $\widetilde w^0=0$).
However, by the cancellation formula~\eqref{cancel}, almost all of the
non-effective weights cancel. Examining the cases, we find that
\begin{align*}
P_{2k+1}(x)&=\sum_j\widetilde w_j^{2k+1}&&N\equiv n\mod 2\\
P_{2k+1}(x)&=\sum_j\widetilde w_j^{2k+1}+(-1)^n\bigl(\tfrac12\bigr)^{2k+1}
&&N\not\equiv n\mod 2.
\end{align*}
If we now apply this formula to the trivial representation, where $N=1$ and
$x=\delta$, we readily obtain the statement of the proposition.
\end{proof}
\begin{cor} For $N$ odd,
\begin{equation} \label{coridodd}
\sum_{j=1}^N \widetilde w_j - \frac{n-1}2 = 0,\qquad
\sum_{j=1}^N (\widetilde w_j)^3 - \Bigl(\frac{n-1}2\Bigr)^3 = 3 c(\lambda),
\end{equation}
and for $N$ even,
\begin{equation} \label{corideven}
\sum_{j=1}^N \widetilde w_j - \frac12 - \frac{n-1}2 = 0,\qquad
\sum_{j=1}^N (\widetilde w_j)^3 - \Bigl(\frac12\Bigr)^3
- \Bigl(\frac{n-1}2\Bigr)^3 = 3 c(\lambda).
\end{equation}\end{cor}

\begin{remk} \label{remparity}
The distinction based on the parity of $N$ (which coincides, for generic
representations, with the parity of the dimension $n$) can be removed by
adding a ``dummy'' conformal weight to the sum: one can either add a
translated conformal weight $\widetilde w=-1/2$ when $N$ is even, or,
following Branson~\cite{branson-elliptic}, a translated conformal weight
$\widetilde w=1/2$ when $N$ is odd. This remark remains true for all the
results proved in this section, provided care is taken in exceptional cases
where the dummy conformal weight already occurs as an effective conformal
weight.
\end{remk}

We now obtain a generating series for the higher Casimirs. These are similar
to the expressions of Perelomov and Popov~\cite{PP}, but differ in two
significant ways: firstly, we compute $\ptr\widetilde B^\ell$, rather than
$\ptr B^\ell$; and secondly, we give the generating series in terms of
translated conformal weights, rather than coordinates of $\lambda$.

\begin{prop}\label{propsup} The partial traces of $\widetilde B^\ell$ are
given by the following generating series:
\begin{equation*}
1 + \sum_{\ell \geq 0} \ptr\widetilde B^{\ell}\, t^{\ell + 1}
= \frac{t}2 + \Bigl( 1 - (-1)^N \frac{t}2\Bigr)
\, \prod_{j = 1}^N \frac{1 + \widetilde w_j t}{1 -  \widetilde w_j t}.
\end{equation*}
\end{prop}
\begin{proof} For each $\ell$,
\begin{equation*}
\ptr\widetilde B^{\ell} = \frac{\trace\widetilde B^\ell}{\dim V}
=\sum (\widetilde w_j)^{\ell}\frac{\dim W_j}{\dim V},
\end{equation*}
since the partial traces act by scalars on $V$. The relative dimensions
$\dim W_j/\dim V$ may be computed as follows.
\begin{lemma} \label{lemmadim}
Let $\Res_{z = \widetilde w_j} (\cdot)$  denote the residue at
$\widetilde w_j$ of the rational function within parentheses. Then:
\begin{enumerate}
\item if $N$ is odd,  
\begin{equation*}
\frac{\dim W_j}{\dim V} = (2 \widetilde w_j + 1) \prod _{k\neq j} 
\frac{\widetilde w_j + \widetilde w_k }{\widetilde w_j - \widetilde w_k} 
= \Res_{z = \widetilde w_j}  
\left( \frac{z + \frac12}{z} \, \prod_{k = 1}^N  
\frac{z + \widetilde w_k }{z  -  \widetilde w_k} \right), 
\end{equation*}
\item if $N$ is even,
\begin{equation*}
\frac{\dim W_j}{\dim V} = (2 \widetilde w_j -  1) \prod_{k \neq j} 
\frac{\widetilde w_j + \widetilde w_k }{\widetilde w_j  -  \widetilde w_k} 
= \Res_{z =  \widetilde w_j}  
\left( \frac{z -  \frac12}{z} \, \prod_{k = 1}^N  
\frac{z + \widetilde w_k }{z  -  \widetilde w_k} \right).
\end{equation*}
\end{enumerate}
\end{lemma}

{\flushleft\it Proof of the lemma.} Weyl's dimension 
formula (see for instance \cite{fulton-harris,samelson}) gives
\begin{equation} \dim W_j = \prod_{\alpha \in \mathcal{R}^+} 
\frac{\langle\mu^{(j)} + \delta, \alpha\rangle}{\langle\delta,
 \alpha\rangle},\qquad  
\dim V = \prod_{\alpha \in \mathcal{R}^+} 
\frac{\langle\lambda + \delta, \alpha\rangle}{\langle\delta, \alpha\rangle},
\end{equation}
where $\mathcal{R}^+$ is the set of positive roots of $\lie{so}(n)$, hence
\begin{equation} \frac{\dim W_j}{\dim V} = 
\prod_{\alpha \in \mathcal{R}^+} \frac{\langle\mu^{(j)} + \delta,
\alpha\rangle}{\langle\lambda + \delta, \alpha\rangle}.
\end{equation}
Unless the dominant weight $\mu^{(j)}$ of $W_j$ is equal to $\lambda$, 
$\mu^{(j)}$ is one of the $2m$ virtual weights $\mu^{i,\pm}=\lambda\pm\eps_i$.
Hence
\begin{equation} \label{dimvirt} 
\frac{\dim W_j}{\dim V} = \prod_{\alpha \in \mathcal{R}^+} 
\Bigl( 1 \pm  \frac{\alpha_i}{\langle\lambda + \delta, \alpha\rangle} \Bigr) 
\end{equation}
so that
\begin{equation}  \label{dimeven}
\frac{\dim W_j}{\dim V} = 
\prod_{\{\widetilde w^{k,\pm}: k \neq i(j)\}}
  \frac{\widetilde w^{i,\pm} +
  \widetilde w^{k,\pm}}{\widetilde w^{i,\pm}  - \widetilde w^{k,\pm}},
\end{equation}
if $n=2m$ is even, and
\begin{equation}  \label{dimodd} 
\frac{\dim W_j}{\dim V}  = \frac{\widetilde w^{i,\pm} +
  \frac12}{\widetilde w^{i,\pm} - \frac12}\, 
  \prod_{\{\widetilde w^{k,\pm}: k \neq i(j)\}} \frac{\widetilde w^{i,\pm}  +
  \widetilde w^{k,\pm}}{\widetilde w^{i,\pm} - \widetilde w^{k,\pm}},
\end{equation}
if $n = 2m + 1$ is odd. Applying the cancellation rule~\eqref{cancel}
and analyzing each case in turn completes the proof.\qed

\smallskip

{\flushleft\it Proof of Proposition
\textup{\ref{propsup} (}continued\textup).}
It follows from the lemma that
\begin{align*}
\frac{\trace \widetilde B^{\ell} }{\dim V}
&= \sum _{j = 1}^N \Res_{z =\widetilde w_j} 
\left( z^{\ell - 1} \Bigl(z - \frac{(-1)^{N}}2\Bigr)
\prod_{k = 1}^N  \frac{z + \widetilde w_k }{z  -  \widetilde w_k} \right)\\
&=\Res_{t=0}\left(t^{-2}t^{1-\ell}\Bigl(\frac1t-\frac{(-1)^{N}}2\Bigr)
\prod_{k = 1}^N  \frac{1/t +  \widetilde w_k }{1/t - \widetilde w_k} \right)
\end{align*}
by the residue theorem. It is straightforward to check that this residue is
the coefficient of $t^{\ell+1}$ in the desired rational expression of
Proposition \ref{propsup}.
\end{proof}
\begin{cor} The partial traces of $\widetilde B^\ell$ are given
by the generating series:
\begin{equation*}
1 + \sum_{\ell \geq 0} \ptr\widetilde B^{\ell}\, t^{\ell + 1}
= \frac{t}2 + \Bigl( 1 - (-1)^N \frac{t}2\Bigr)S(t)
\end{equation*}
where $S'(t)/S(t) = 2 \sum s_{2k+1}(\widetilde w) t^{2k+1}$ and
$s_{2k+1}(\widetilde w)$ are the power sum symmetric functions in the
translated conformal weights.  In particular, by
Proposition~\textup{\ref{propid}}, the partial traces can be computed from the
polynomials $P_{2k+1}(x)$.
\end{cor}
We recover from these generating functions, the results of Perelomov
and Popov for the orthogonal Lie algebras~\cite{PP}. Although the generating
functions are not too complicated, they suggest that the operators
$\widetilde A_k$ defined by
\begin{equation}\widetilde A_k
=\sum_{\ell=0}^k(-1)^\ell\sigma_\ell(\widetilde w)\widetilde B^{k-\ell},
\end{equation}
where $\sigma_\ell(\widetilde w)$ denotes the $\ell$th elementary symmetric
function in the translated conformal weights, will have much simpler
traces. This is indeed the case.

\begin{prop} The partial trace of $\widetilde A_j$ is:
\begin{equation}
\ptr\widetilde A_j=\bigl(1+(-1)^j)\sigma_{j+1}(\widetilde w)+
\frac12\bigl((-1)^j-(-1)^N\bigr)\sigma_j(\widetilde w).
\end{equation}
\proofof{prop} We compute the generating function
\begin{align*}
\sum_{j\geq0} \ptr\widetilde A_j t^{j+1}&=\sum_{j\geq0}\sum_{k=0}^j(-1)^k
\sigma_k(\widetilde w)\ptr\widetilde B^{j-k} t^{j+1}\\
&=\sum_{k\geq0}\sum_{j\geq k}(-1)^k\sigma_k(\widetilde w)
\ptr\widetilde B^{j-k} t^kt^{j-k+1}\\
&=\sum_{k\geq0}(-1)^k\sigma_k(\widetilde w)t^k\sum_{\ell\geq0}
\ptr\widetilde B^{\ell} t^{\ell+1}\\
&=\Bigl(1-\frac{(-1)^Nt}2\Bigr)\prod_{j=1}^N(1+\widetilde w_jt)
-\Bigl(1-\frac t2\Bigr)\prod_{j=1}^N(1-\widetilde w_jt)\\
&=\sum_{j\geq0} \Bigl( \bigl(1-(-1)^j\bigr)
+\frac t2\bigl((-1)^j-(-1)^N\bigr)\Bigr)\sigma_j(\widetilde w)t^j.
\end{align*}
This yields the stated formula.
\end{proof}

We are now ready for the main result of this section.

\begin{thm}\label{maincas} Define
$\widetilde C_j=\widetilde A_j+\frac14\bigl((-1)^N-(-1)^j\bigr)
\widetilde A_{j-1}$, where $\widetilde A_{-1}=0$ by convention. Then
$(\widetilde C_j)_{\alpha\tens\beta}
=(-1)^j(\widetilde C_j)_{\beta\tens\alpha}$.
\end{thm}
\begin{cor}\label{maincor} If $N$ is odd then
\begin{equation}
\ip{\widetilde A_{2j+1}(\alpha\tens v),\alpha\tens v}=0
\end{equation}
while if $N$ is even,
\begin{equation}
\ip{\widetilde A_{2j+1}(\alpha\tens v),\alpha\tens v}+
\frac12\ip{\widetilde A_{2j}(\alpha\tens v),\alpha\tens v}=0.
\end{equation}
\end{cor}

The idea of looking for polynomials in $B$ with symmetry properties
was first suggested to the authors by T. Diemer and G. Weingart
(private communication). One of their key results is the following:

\begin{thm}\textup{(Diemer-Weingart)}\label{Tammo-Gregor}
Let $q_j(B)$ be a sequence of polynomials in $B$ with $q_j(B)=0$ for
$j<0$, $q_0(B)=1$ and for $j\geq0$,
\begin{equation}\begin{split}
q_{j+1}(B)_{\alpha\tens\beta}
&=\biggl(\Bigl(B+\frac{n-1+(-1)^j}2\iden\Bigr)\circ q_j(B)\biggl)
_{\alpha\tens\beta}\\
&\qquad-\frac12\ip{\alpha,\beta}\ptr q_j(B)
+\sum_{k\geq1}a_{jk}\,q_{j+1-2k}(B)_{\alpha\tens\beta}
\end{split}\end{equation}
for some $a_{jk}\in\R$. Then
\begin{equation}\label{concl}
q_j(B)_{\alpha\tens\beta}=(-1)^jq_j(B)_{\beta\tens\alpha}.
\end{equation}
\proofof{thm} We give the proof of Diemer and
Weingart, which is by complete induction on $j$:
clearly~\eqref{concl} holds for $j\leq0$ and we have an inductive
formula for $q_{j+1}$. Introducing the temporary notation
$(c_j)_{\alpha\tens\beta}=\frac12\ip{\alpha,\beta}\bigl(\ptr q_j(B)\bigr)$
we have
\begin{align*}
&2\bigl(q_{j+1}(B)_{\alpha\tens\beta}
-(-1)^{j+1}q_{j+1}(B)_{\beta\tens\alpha}\bigr)\\
&\quad=\Bigl(\bigl(2 B+(n-1+(-1)^j)\iden\bigr)\circ q_j(B)
-c_j\Bigr)_{\alpha\tens\beta}\\
&\qquad+(-1)^j\Bigl(\bigl(2B+(n-1+(-1)^j)\iden\bigr)\circ q_j(B)
-c_j\Bigr)_{\beta\tens\alpha}\\
&\quad=(B\circ q_j(B))_{\alpha\tens\beta}
+(-1)^j(q_j(B)\circ B)_{\beta\tens\alpha}\\
&\qquad+(-1)^j\bigl((B\circ q_j(B))_{\beta\tens\alpha}
+(-1)^j(q_j(B)\circ B)_{\alpha\tens\beta}\bigr)\\
&\qquad+\bigl((n-1+(-1)^j)q_j(B)-c_j\bigr)_{\alpha\tens\beta}
+(-1)^j\bigl((n-1+(-1)^j)q_j(B)-c_j\bigr)_{\beta\tens\alpha}
\end{align*}
since $q_j(B)$ commutes with $B$. The result follows by observing that
\begin{align*}
(B\circ q_j(B))_{\alpha\tens\beta}&+(-1)^j(q_j(B)\circ B)_{\beta\tens\alpha}\\
&=\sum_i\bigl(B_{\alpha\tens e_i}\circ q_j(B)_{e_i\tens\beta}
+(-1)^j q_j(B)_{\beta\tens e_i}\circ B_{e_i\tens\alpha}\bigr)\\
&=\sum_i\bigl(B_{\alpha\tens e_i}\circ q_j(B)_{e_i\tens\beta}
-q_j(B)_{e_i\tens\beta}\circ B_{\alpha\tens e_i} \bigr)\\
&=\sum_i[d\lambda(\alpha\wedge e_i),q_j(B)_{e_i\tens\beta}]
=\sum_iq_j(B)_{\alpha\wedge e_i.(e_i\tens\beta)}
\end{align*}
by equivariance of $q_j(B)$, where $\alpha\wedge e_i.(e_i\tens\beta)$ is
defined using the action of $\lie{so}(n)$ on $\R^n\tens\R^n$. This gives,
finally,
\begin{align*}
(B\circ q_j(B))_{\alpha\tens\beta}&+(-1)^j(q_j(B)\circ B)_{\beta\tens\alpha}\\
&=q_j(B)_{\alpha\tens\beta}-n q_j(B)_{\alpha\tens\beta}
+\ip{\alpha,\beta}\ptr q_j(B)-q_j(B)_{\beta\tens\alpha}\\
&=\bigl((1-n-(-1)^j)q_j(B)+c_j\bigr)_{\alpha\tens\beta},
\end{align*}
which completes the proof.
\end{proof}

By taking $a_{jk}=0$ (for all $j,k$), Diemer and Weingart obtain an inductive
definition of a sequence of polynomials with the desired symmetry properties.
Unfortunately, the task of computing these polynomials explicitly is
formidable because of the complexity of the traces of the powers of $B$.

The polynomials $\widetilde C_j$ defined here are completely explicit and
because they have simple traces we are able to prove that they satisfy the
inductive conditions of Theorem~\ref{Tammo-Gregor}.  More precisely, we have:
\begin{lemma} For $j\geq0$,
\begin{align*} \widetilde C_{j+1}
&=\Bigl(\widetilde B+\frac{(-1)^j}2\iden\Bigr)\circ\widetilde C_j
-\frac12\ptr\widetilde C_j\\
&\quad+\tfrac18\bigl(1-(-1)^{N+j}\bigr)\widetilde C_{j-1}
+\tfrac12\bigl(1-(-1)^j\bigr)\Bigr(\sigma_{j+1}(\widetilde w)
-\tfrac12\bigl(1-(-1)^N\bigr)\sigma_j(\widetilde w)\Bigr)\iden
\end{align*}
\proofof{lemma} Note that
$\widetilde C_j=\widetilde A_j+\frac14\bigl((-1)^N-(-1)^j\bigr)
\widetilde C_{j-1}$ and so
\begin{align*}
&\widetilde C_{j+1}-\widetilde B\widetilde C_j-\tfrac12(-1)^j\widetilde C_j
=\widetilde A_{j+1}-\widetilde B\widetilde A_j-\tfrac12(-1)^j\widetilde C_j
+\tfrac14\bigl((-1)^N+(-1)^j\bigr)\widetilde C_j\\
&\hspace{9cm}-\tfrac14\bigl((-1)^N-(-1)^j\bigr)\widetilde B\widetilde C_{j-1}\\
&\quad=\widetilde A_{j+1}-\widetilde B\widetilde A_j
+\tfrac14\bigl((-1)^N-(-1)^j\bigr)
\bigl(\widetilde C_j-\widetilde B\widetilde C_{j-1}\bigr)\\
&\quad=\widetilde A_{j+1}-\widetilde B\widetilde A_j
+\tfrac14\bigl((-1)^N-(-1)^j\bigr)
\bigl(\widetilde C_j-\widetilde B\widetilde C_{j-1}
-\tfrac12(-1)^j\widetilde C_{j-1}\bigr)\\
&\hspace{9cm}+\tfrac18\bigl(1-(-1)^{N+j}\bigr)\widetilde C_{j-1}\\
&\quad=\widetilde A_{j+1}-\widetilde B\widetilde A_j
+\tfrac14\bigl((-1)^N-(-1)^j\bigr)
\bigl(\widetilde A_j-\widetilde B\widetilde A_{j-1}\bigr)
+\tfrac18\bigl(1-(-1)^{N+j}\bigr)\widetilde C_{j-1}.
\end{align*}
Now, by definition, we have
$\widetilde A_{j+1}-\widetilde B\widetilde A_j
=(-1)^{j+1}\sigma_{j+1}(\widetilde w)\iden$
and so
\begin{equation}\label{step}\begin{split}
\widetilde C_{j+1}-\widetilde B\widetilde C_j-\tfrac12(-1)^j\widetilde C_j
&=\tfrac18\bigl(1-(-1)^{N+j}\bigr)\widetilde C_{j-1}\\
&\quad-(-1)^j\sigma_{j+1}(\widetilde w)\iden-
\tfrac14\bigl(1-(-1)^{N+j}\bigr)\sigma_j(\widetilde w)\iden.
\end{split}\end{equation}
Finally, observe that
\begin{equation*}
\ptr\widetilde C_j=\bigl(1+(-1)^j\bigr)\bigl(\sigma_{j+1}(\widetilde w)
+\tfrac14(1-(-1)^{N+j})\sigma_j(\widetilde w)\bigr)\iden.
\end{equation*}
Adding one half of this onto~\eqref{step} completes the proof.
\end{proof}
Theorem~\ref{maincas} follows immediately from this Lemma
and Theorem~\ref{Tammo-Gregor}.

\section{Refined Kato inequalities}\label{mainsect}

In the last section we learnt that by working with $\widetilde B$ and
$\widetilde A_j$ instead of $B$ and $A_j$, we could obtain some explicit
formulae. Of course $\widetilde B=B+\frac12(n-1)\iden\,$ has the same
eigenspaces as $B$ and so we can rewrite~\eqref{proj} as:
\begin{equation}\label{projt}
|\Pi_j(\alpha\tens v)|^2=\frac{\displaystyle\sum_{k=0}^{N-1}\widetilde
w_j^{N-1-k}
\ip{\widetilde A_k(\alpha\tens v),\alpha\tens v}}
{\displaystyle\prod_{k\neq j}(\widetilde w_j-\widetilde w_k)}.
\end{equation}

If $N$ is odd, Corollary~\ref{maincor} implies that the terms with $k$ odd
vanish, while for $N$ even, we have
\begin{equation*}
\ip{\widetilde A_{2j+1}(\alpha\tens v),\alpha\tens v}+
\frac12\ip{\widetilde A_{2j}(\alpha\tens v),\alpha\tens v}=0.
\end{equation*}
Our main result will readily follow from this.

\begin{maintheo}\label{thmain}
Let $I$ a subset of $\{1,\ldots,N\}$ corresponding to an operator $P_I$ acting
on $E$. Then a Kato constant $k_I$ for the kernel of $P_I$ is given by the
following expressions.

If $N$ is odd, then
\begin{equation}\label{Nodd} k_I^2 = \max_{J\in\NE}
\Biggl(\,\sum_{i\in\widehat I\intersect\widehat J} \frac{\prod_{j\in J} (\widetilde w_i
+ \widetilde w_j)}{\prod_{j\in\widehat J\setminus\{i\}} (\widetilde w_i -
\widetilde w_j)} \Biggr)
= 1 - \min_{J\in\NE} \Biggl(\,\sum_{i\in I\intersect\widehat J} \frac{\prod_{j\in
J} (\widetilde w_i + \widetilde w_j)}{\prod_{j\in \widehat J\setminus\{i\}}
(\widetilde w_i - \widetilde w_j)} \Biggr).
\end{equation}

If $N$ is even, then
\begin{multline}\label{Neven} k_I^2 = \max_{J\in\NE}
\Biggl(\,\sum_{i\in\widehat I\intersect\widehat J} \frac{ \bigl(\widetilde
w_i-\tfrac12\bigr) \prod_{j\in J} (\widetilde w_i + \widetilde
w_j)}{\prod_{j\in \widehat J\setminus\{i\}} (\widetilde w_i - \widetilde w_j)}
\Biggr)\\ = 1 - \min_{J\in\NE} \Biggl(\,\sum_{i\in I\intersect\widehat J}
\frac{\bigl(\widetilde w_i-\tfrac12\bigr)\prod_{j\in J} (\widetilde w_i +
\widetilde w_j)}{\prod_{j\in \widehat J\setminus\{i\}} (\widetilde w_i -
\widetilde w_j)} \Biggr).
\end{multline}
These constants are sharp, unless $N=2\nu+1$, $\lambda$ is properly
half-integral, and the set $J$ achieving the extremum contains $\nu+1$.
\end{maintheo}
Recall that $\NE$ denotes the set of subsets of $\{1,\ldots N\}$ whose
elements are obtained by choosing exactly one index in each of the sets
$\{j,N+2-j\}$ for each $j$ with $2\leq j\leq\nu$ if $N=2\nu-1, 2\nu$ and for
each $j$ with $2\leq j\leq\nu+1$ if $N=2\nu+1$. These correspond to the
maximal non-elliptic operators unless $N=2\nu+1$ and $\lambda$ is properly
half-integral, when there are also some elliptic subsets in $\NE$.

Explicit values of the constants for a number of cases, including all minimal
elliptic operators, will be given in sections 6 and 7, and in the
appendix. Note that $k_I=1$ for non-elliptic operators, as one would expect.

{\flushleft\bf Proof of the Main Theorem}.
We let first $N=2\nu-1$ and denote
$Q_k=(-1)^{k-1}\ip{\widetilde A_{2k-2}(\alpha\tens v),\alpha\tens v}$. We have
\begin{equation}\label{projodd}
|\Pi_j(\alpha\tens v)|^2=\frac{\displaystyle\sum_{k=1}^{\nu}
\widetilde w_j^{2(\nu-k)}(-1)^{k-1}Q_k}
{\displaystyle\prod_{k\neq j}(\widetilde w_j-\widetilde w_k)}
=\frac{\displaystyle\widetilde w_j^{2(\nu-1)}-
\sum_{k=2}^{\nu}\widetilde w_j^{2(\nu-k)}(-1)^kQ_k}
{\displaystyle\prod_{k\neq j}(\widetilde w_j-\widetilde w_k)}
\end{equation}
since $Q_1=1$. We can now obtain bounds on $Q_2,\ldots Q_\nu$
using the non-negativity of the norms. Since the denominator
in ~\eqref{projodd} has sign $(-1)^{j-1}$ these inequalities are:
\begin{equation}\label{oddsys}
\sum_{k=2}^{\nu}(-1)^{j+k}\widetilde w_j^{2(\nu-k)}Q_k
\geq(-1)^j\widetilde w_j^{2(\nu-1)}
\end{equation}
with equality iff $|\Pi_j(\alpha\tens v)|^2=0$. 

This system of linear inequalities confines the values of the $Q_k$'s to a
convex region in $\R^{\nu-1}$. Our first goal is to show that this region is
compact, hence polyhedral, and to identify its vertices. For this we let
$\pi_j$ denote the affine functions of $Q=(Q_2,\ldots Q_\nu)$ given by
$|\Pi_j(\alpha\tens v)|^2$ and note the following.
\begin{lemma}\label{valuespiq} Let $J$ be a subset of $\{1,...,N\}$ with
$\nu-1$ elements. Then the intersection of the $\nu-1$ affine hyperplanes
$\pi_j=0$ for all $j\in J$ consists of the single point $Q_J=(Q_2,\ldots
Q_\nu)$ with $Q_k = \sigma_{k-1}\bigl((\widetilde{w}_j^2)_{j\in J}\bigr)$.
At this point the affine functions $\pi_j$ take the values
\begin{equation}\label{oddsoln}
\pi_j(Q_J)=
\frac{\displaystyle\prod_{k\in J}(\widetilde w_j^2-\widetilde w_k^2)}
{\displaystyle\prod_{k\neq j}(\widetilde w_j-\widetilde w_k)}
=\frac{\displaystyle\prod_{k\in J, k\neq j}(\widetilde w_j+\widetilde w_k)}
{\displaystyle\prod_{k\in \widehat J, k\neq j}(\widetilde w_j-\widetilde w_k)}
\,\varepsilon_j(J)
\end{equation}
where $\varepsilon_j(J)=0$ if $j\in J$ and $1$ if not.
\end{lemma}
This lemma follows simply by observing that the affine function $\pi_j$ is
obtained by evaluating a polynomial independent of $j$ on $\widetilde w_j^2$,
and then using the fact that the coefficients of a polynomial are the
elementary symmetric functions of the roots.

Compactness of the convex region is now obtained by taking
$J=\{2,\ldots\nu\}$ and $J=\{\nu+1,\ldots 2\nu-1\}$. The inverse of the
Vandermonde system of inequalities for $J=\{2,\ldots\nu\}$ has non-negative
entries, while for $J=\{\nu+1,\ldots 2\nu-1\}$, it has non-positive
entries.
\begin{prop} Let $N=2\nu-1$. Then for $k=2,\ldots\nu$,
\begin{equation}\label{oddbound}
\sigma_{k-1}(\widetilde w_2^2,\ldots\widetilde w_\nu^2)\leq Q_k
\leq\sigma_{k-1}(\widetilde w_{\nu+1}^2,\ldots \widetilde w_{2\nu-1}^2).
\end{equation}
The lower bounds are all attained if and only if
$\Pi_{\{2,\ldots\nu\}}(\alpha\tens v)=0$, while the upper bounds are all
attained if and only if $\Pi_{\{\nu+1,\ldots2\nu-1\}}(\alpha\tens v)=0$.
These bounds are sharp by non-ellipticity of $P_{\{2,\ldots\nu\}}$ and
$P_{\{\nu+1,\ldots2\nu-1\}}$.
\end{prop}
When $N=2\nu-1=3$, the case most commonly occuring in practice, it is now
straightforward to obtain sharp Kato constants. However, for $N\geq5$, the
upper bound for some $Q_k$ and the lower bound for another (as given in this
proposition) will not be simultaneously attained: the convex region is
smaller. We illustrate this in the case $N=5$ ($\nu=3$).

In this diagram, the numbered lines represent the conditions on $Q_2$ and
$Q_3$ for the norms of $\Pi_1,\ldots \Pi_5$ to vanish. The shaded region
represents the range of possible values for $(Q_2,Q_3)$, while the dotted
rectangle represents the bounds on $(Q_2,Q_3)$ we have found. We have circled
the points corresponding to the non-elementary minimal elliptic operators.

\medbreak \input{epsf} \centerline{\epsfbox{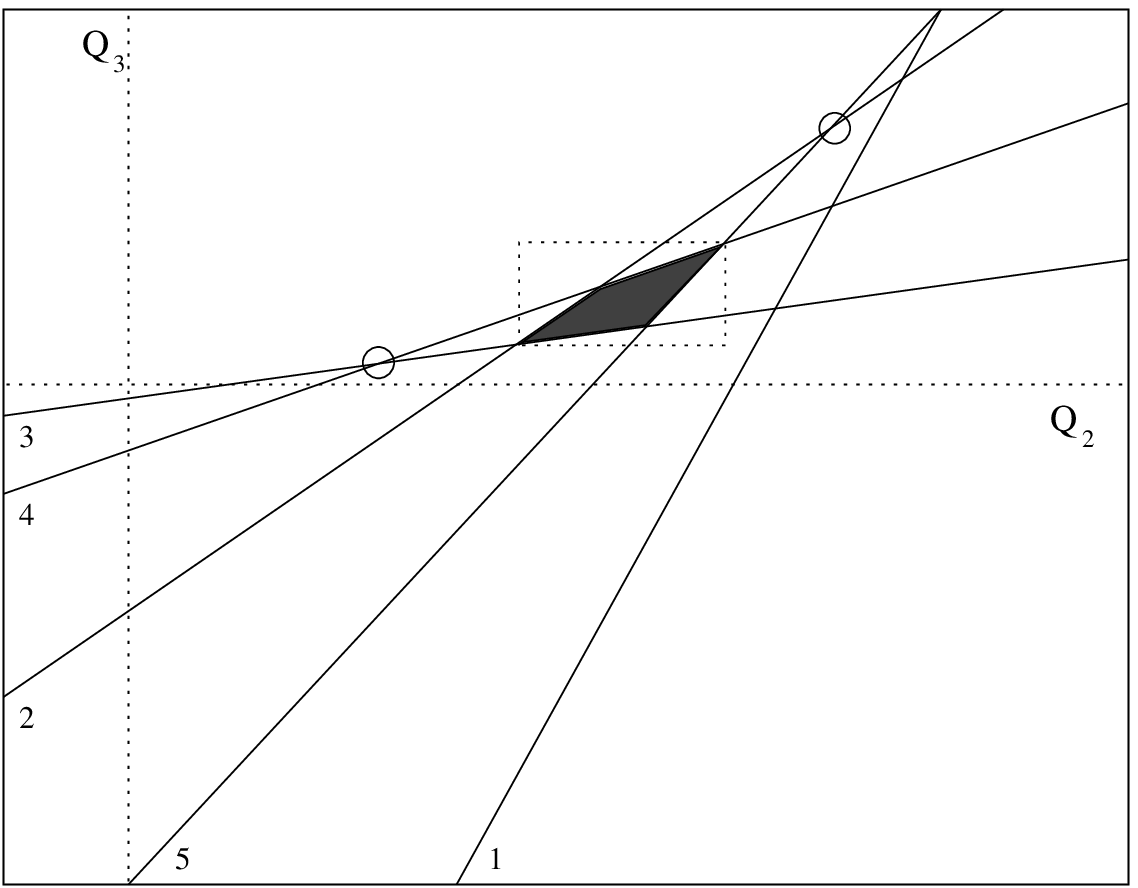}} \medbreak

According to Ansatz~\ref{ansatz}, in order to find a sharp Kato constant for
$P_I$ we must maximize (for $|\alpha\tens v|=1$) the projection
$|\Pi_{\widehat I}(\alpha\tens v)|^2=1-|\Pi_I(\alpha\tens v)|^2$, which is
equivalent to minimizing $|\Pi_I(\alpha\tens v)|^2=\sum_{i\in I}\pi_i$.

Since these norms are affine in the $Q_k$'s, it follows that to minimize or
maximize them on the polyhedral region of admissible values of the $Q_k$'s,
we must find the supporting hyperplanes associated to the linear part of the
function. Such a supporting hyperplane certainly contains a vertex of the
polyhedron, and so it suffices to minimize or maximize over the set of
vertices.

We claim that these vertices are the points $Q_J$ with $J\in\NE$.  Certainly
these points are vertices, since if $J\in\NE$ then $P_J$ is non-elliptic (this
part of the argument will fail when $N=2\nu+1$) and so there is some
$\alpha\tens v$ of norm one with $\Pi_j(\alpha\otimes v) = 0$ for each $j$ in
$J$. Therefore it remains to eliminate the points $Q_J$ with $J\notin\NE$ as
possible vertices, which we do by showing that a point $Q_J$ with $P_J$
elliptic does not lie in the polyhedral region. This is done by proving that
there is, for every such $J$, an index $i$ such that the affine function
$\pi_i$ assumes a (strictly) negative value at $Q_J$.
Equation~\eqref{oddsoln} tells us that for $i\notin J$, $\pi_i(Q_J)$ is
nonzero and has the sign $(-1)^{i-1}\rho_i$ where $\rho_i$ is the sign of
$\prod_{j\in J}( \widetilde{w}_i^2 - \widetilde{w}_j^2 )$.  If $P_J$ is
elliptic, $J$ contains a minimal elliptic set, hence either the index $1$ or a
couple of indices of the form $(j,N+2-j)$. In any case, since $J$ has length
$\nu-1$ and there are exactly $\nu-1$ couples of type $(\ell,N+2-\ell)$, there
is at least one such couple \emph{outside} $J$.  One readily checks that
$\widetilde w_{\ell}^2$ and $\widetilde w_{N+2-\ell}^2$ are adjacent in the
ordering of the squares of the conformal weights, and so
$\rho_{\ell}=\rho_{N+2-\ell}$. Since $N$ is odd, $\ell$ and $N+2-\ell$ have
the opposite parity, and so one of $i=\ell$ or $i=N+2-\ell$ yields a negative
sign for $\pi_i$. This proves the claim, and now maximizing or minimizing over
the vertices using~\eqref{oddsoln} proves the main theorem for $N=2\nu-1$.

The argument for the case $N=2\nu+1$ is completely analogous, by replacing
$\nu$ with $\nu-1$.  When $\lambda$ properly half-integral, the lower bounds
in the analogue of~\eqref{oddbound} will not be sharp since
$P_{\{2,\ldots\nu+1\}}$ is elliptic. However, we only used these bounds to
establish compactness of the convex region defined by the nonnegativity of the
norms, so this does not matter. The ellipticity of $P_{\nu+1}$ means that some
of the vertices of this polyhedral region are not possible values for the
$Q_k$'s. More precisely, the index sets corresponding to the vertices are
still contained in the set $\NE$, and so we can maximize or minimize over
$\NE$, but we will not obtain sharp results if the extremum is obtained
at a vertex corresponding to an index set containing $\nu+1$.

Now suppose $N=2\nu$ and let $Q_k=(-1)^{k-1}\ip{\widetilde
A_{2k-2}(\alpha\tens v),\alpha\tens v}$ we have
\begin{equation}\label{projeven}\begin{split}
|\Pi_j(\alpha\tens v)|^2&=\frac{\widetilde w_j-\frac12}
{\displaystyle\prod_{k\neq j}(\widetilde w_j-\widetilde w_k)}
\sum_{k=1}^{\nu}\widetilde w_j^{2(\nu-k)}(-1)^{k-1}Q_k\\
&=\frac{\widetilde w_j-\frac12}
{\displaystyle\prod_{k\neq j}(\widetilde w_j-\widetilde w_k)}
\biggl(\widetilde w_j^{2(\nu-1)}-
\sum_{k=2}^{\nu}\widetilde w_j^{2(\nu-k)}(-1)^kQ_k\biggr)
\end{split}\end{equation}
since $Q_1=1$. Our strategy is now the same as before: we obtain the
polyhedron using the non-negativity of the norms and its vertices by looking
at maximal length non-elliptic operators.
Since the denominator in~\eqref{projeven} has sign
$(-1)^{j-1}$ these inequalities are:
\begin{equation}\label{evensys}\begin{split}
\sum_{k=2}^{\nu}(-1)^{j+k}\widetilde w_j^{2(\nu-k)}Q_k
&\geq(-1)^j\widetilde w_j^{2(\nu-1)}\qquad\mathrm{for}\quad j\leq\nu\\
\sum_{k=2}^{\nu}(-1)^{j+k}\widetilde w_j^{2(\nu-k)}Q_k
&\leq(-1)^j\widetilde w_j^{2(\nu-1)}\qquad\mathrm{for}\quad j\geq\nu+1.
\end{split}\end{equation}
Lemma~\ref{valuespiq} is unchanged except that the formula for $\pi_i(Q_J)$
have an additional $\widetilde w_i-\frac12$. To obtain compactness, we
consider $J=\{2,\ldots\nu\}$ and $J=\{\nu+2,\ldots 2\nu\}$ and again observe
that the inverses of these Vandermonde systems have entries all of one sign.
\begin{prop} Let $N=2\nu$. Then for $k=2,\ldots\nu$,
\begin{equation}\label{evenbound}
\sigma_{k-1}(\widetilde w_2^2,\ldots\widetilde w_\nu^2)\leq Q_k
\leq\sigma_{k-1}(\widetilde w_{\nu+2}^2,\ldots \widetilde w_{2\nu}^2).
\end{equation}
The lower bounds are all attained if and only if
$\Pi_{\{2,\ldots\nu\}}(\alpha\tens v)=0$, while
the upper bounds are all attained if and only if
$\Pi_{\{\nu+2,\ldots2\nu\}}(\alpha\tens v)=0$.
These bounds are sharp by non-ellipticity of
$P_{\{2,\ldots\nu\}}$ and $P_{\{\nu+2,\ldots2\nu\}}$.
\end{prop}
The vertices are identified with $\NE$ in a similar way to the case
$N=2\nu-1$. The only difference comes from the way sign changes when passing
from $i=ell$ to $i=N+2-\ell$: the parity of $i$ does not change but the sign
of the factor $\widetilde{w}_i-1/2$ does.  This proves the main theorem for
$N=2\nu$. \qed

\medbreak

In the next two sections we shall calculate some of the constants more
explicitly, by finding the vertex at which the maximum or minimum is
achieved. This is only feasible when the number of terms in the sum is small
and in general, the vertex depends on the coordinates of $\lambda$.
Nevertheless, this is a worthwhile task, as explicit constants are of more
practical use than extrema over exponentially large sets.

Our main tool is the order of the conformal weights, together with the fact
that, for $j\in\{2,\ldots\nu\}$, we have $\widetilde w_j+\widetilde w_{N+2-j}
=k_j-k_{j-1}<0$. Similarly, for $N=2\nu+1$,
$\widetilde w_{\nu+1}+\widetilde w_{\nu+2}=\widetilde w_{\nu+2}=-\lambda_m<0$.
Hence for any $i\in\{1,\ldots N\}$ and $j\in\{2\ldots\nu\}$:
\begin{multline*}
(\widetilde w_i+\widetilde w_j)(\widetilde w_i-\widetilde w_j)
-(\widetilde w_i+\widetilde w_{N+2-j})(\widetilde w_i-\widetilde w_{N+2-j})\\
=-(\widetilde w_j+\widetilde w_{N+2-j})(\widetilde w_j-\widetilde w_{N+2-j})>0
\end{multline*}
and this also holds for $N=2\nu+1$ and $j=\nu+1$.

By considering the possible signs of the terms, we obtain:
\begin{prop}\label{useful} For any $i\in\{1,\ldots N\}$ and
$j\in\{2\ldots\nu\}$ \textup(or $j=\nu+1$ when $N=2\nu+1$\textup)
with $i\neq j$ and $i\neq N+2-j$, we have:
\begin{align*}
\frac{\widetilde w_i+\widetilde w_j}{\widetilde w_i-\widetilde w_{N+2-j}}
>\frac{\widetilde w_i+\widetilde w_{N+2-j}}{\widetilde w_i-\widetilde w_j}>0
&\qquad\textup{iff}\quad i<j\quad\textup{or}\quad N+2-j<i\\
\frac{\widetilde w_i+\widetilde w_{N+2-j}}{\widetilde w_i-\widetilde w_j}
>\frac{\widetilde w_i+\widetilde w_j}{\widetilde w_i-\widetilde w_{N+2-j}}>0
&\qquad\textup{iff}\qquad\;\; j<i<N+2-j.
\end{align*}
\end{prop}

\section{Refined Kato inequalities with $N$ odd}

When $N$ is odd, we have to minimize or maximize over $J\in\NE$, a sum of a
subset of the following terms:
\begin{align*}
\frac{\prod_{j\in J} (\widetilde w_i+\widetilde w_j)}
{\prod_{j\in\widehat J\setminus\{i\}} (\widetilde w_i-\widetilde w_j)}
&=\frac{\widetilde w_i+\widetilde w_{N+2-i}}{\widetilde w_i-\widetilde w_1}
\prod_{\substack{j\in J\\ j\neq N+2-i}}
\frac{\widetilde w_i+\widetilde w_j}{\widetilde w_i -\widetilde w_{N+2-j}}
&&\qquad\textrm{for}\quad i\in \widehat J\setminus\{1\}\\
\frac{\prod_{j\in J} (\widetilde w_1+\widetilde w_j)}
{\prod_{j\in\widehat J\setminus\{i\}} (\widetilde w_1-\widetilde w_j)}
&=\prod_{j\in J}\frac{\widetilde w_1+\widetilde w_j}
{\widetilde w_1-\widetilde w_{N+2-j}}&&
\end{align*}
Using Proposition~\ref{useful}, the first expression is minimized (subject to
$J\not\ni i$) by $J_i^{\min}=\{2,\ldots i-1,N+2-\nu,\ldots N+2-i\}$ (together
with $\nu+2$ if $N=2\nu+1$) and is maximized by
$J_i^{\max}=\{i+1,\ldots\nu,N+2-i,\ldots N\}$ (together with $\nu+1$ if
$N=2\nu+1$).  The second expression is minimized by
$J_1^{\min}=\{N+2-\nu,\ldots N\}$ (together with $\nu+2$ if $N=2\nu+1$) and
maximized by $J_1^{\max}=\{2,\ldots\nu\}$ (together with $\nu+1$ if
$N=2\nu+1$).

This information suffices to find Kato constants for the elementary elliptic
operators and the complements of generalized gradients. Note that
$J_i^{\min}=J_{N-1-i}^{\min}$ and $J_i^{\max}=J_{N-1-i}^{\max}$, which will
give a few more explicit results.

We shall now show how the values of the constants can be computed for the
non-elementary (\ie, length $2$) minimal elliptic operators.

Let $I=\{i,N+2-i\}$ for $i\in\{2,\cdots,\nu\}$ (or $i=\nu+1$ when
$N=2\nu+1$). Then for any $J\in\NE$, $J\intersect I$ has precisely one
element, and hence so does $J\intersect\widehat I$. Therefore, for each $J$,
the sum has only one term, indexed by either $i$ or $N-2-i$, and so the
minimum, over all $J$, is given by the minimum over $J_i^{\min}$ and
$J_{N+2-i}^{\min}$. Unfortunately, each of these two quantities may be the
smallest, depending on the precise values of the conformal weights, so that we
are forced to keep an minimum in our formulas. However, if $N=2\nu-1$ and
$i=\nu$, then the following argument, together with the fact that $\widetilde
w_1-\widetilde w_{\nu+1}>\widetilde w_1-\widetilde w_\nu$, shows that the
minimum is obtained by using $J_{\nu+1}^{\min}$.
\begin{lemma}\label{ratio}  For each $k=1,\ldots\nu-2$
\begin{equation*}
(\widetilde w_\nu+\widetilde w_{k+1})(\widetilde w_{\nu+1}-\widetilde
w_{2\nu-k}) >(\widetilde w_{\nu+1}+\widetilde
w_{k+1})(\widetilde w_\nu-\widetilde w_{2\nu-k}) >0.
\end{equation*}
\proofof{lemma} Positivity holds because $2\nu-k>\nu+1$, while the inequality
follows from the identity
\begin{multline*}
(\widetilde w_\nu+\widetilde w_{k+1})(\widetilde w_{\nu+1}-\widetilde
w_{2\nu-k}) -(\widetilde w_{\nu+1}+\widetilde w_{k+1})(\widetilde
w_\nu-\widetilde w_{2\nu-k})\\ =-(\widetilde w_{k+1}+\widetilde
w_{2\nu-k})(\widetilde w_\nu-\widetilde w_{\nu+1})
\end{multline*}
and the fact that $\widetilde w_{k+1}+\widetilde w_{2\nu-k}<0$.
\end{proof}
A similar argument works when $N=2\nu+1$ and $i=\nu+1$.

We summarize these observations in the following results.
\begin{thm} Let $E$ be associated to a representation $\lambda$ with
$N=2\nu-1$ and let $P_I$ an elliptic operator on sections of $E$ associated to
a subset $I$ of $\{1,\ldots N\}$. Then in the following cases, a refined Kato
inequality of the type $\bigl|d|\xi|\bigr| \leq k_I |\nabla\xi|$ holds outside
the zero set of $\xi$ for $\xi$ in the kernel of $P_I$.
\begin{enumerate}
\item
For $\{1\} \subseteq I \subseteq \{1,\nu+1,\ldots 2\nu-1\}$, we have
\begin{equation*}
k_I^2 = 1 - \frac{\prod_{k=\nu+1}^{2\nu-1}(\widetilde w_1+\widetilde w_k)}
{\prod_{k=2}^{\nu}(\widetilde w_1-\widetilde w_k)}.
\end{equation*}
Equality holds iff $\nabla\xi=\Pi_{\{2,...\nu\}}(\alpha\tens\xi)$ for a
$1$-form $\alpha$ with $\Pi_{\{\nu+1,...2\nu-1\}}(\alpha\tens\xi)=0$.
\item
For $\{i,2\nu+1-i\} \subseteq I \subseteq \{i,2\nu+1-i\}\cup J_0$, with $i\in
\{2,\ldots,\nu\}$ and $J_0=\{j: 2\leq j < i\}\cup\{2\nu+1-j:\ i<j\leq\nu\}$,
we have
\begin{equation*} k_I^2 = 1 - \min (C_1,C_2), \end{equation*}
where
\begin{equation*}\begin{split}
C_1 &= \frac{\widetilde w_i+\widetilde w_{2\nu+1-i}}{\widetilde w_i
-\widetilde w_1}\prod_{k\in J_0}\frac{\widetilde w_i+\widetilde w_k}
{\widetilde w_i-\widetilde w_{2\nu+1-k}}, \\ 
C_2 &= \frac{\widetilde w_i+\widetilde w_{2\nu+1-i}}{\widetilde w_{2\nu+1-i}
-\widetilde w_1}\prod_{k\in J_0}\frac{\widetilde w_{2\nu+1-i}+\widetilde w_k}
{\widetilde w_{2\nu+1-i}-\widetilde w_{2\nu+1-k}}.
\end{split}\end{equation*}
Equality holds iff
$\nabla\xi=\Pi_{\widehat J^0\setminus\{i,2\nu+1-i\}}(\alpha\tens\xi)$ for a
$1$-form $\alpha$ with
\begin{equation*}
\Pi_{\{i\}\cup J_0}(\alpha\tens\xi)=0\,\textrm{ if }\, C_2<C_1\quad
\textrm{or}\quad
\Pi_{\{2\nu+1-i\}\cup J_0}(\alpha\tens\xi)=0\,\textrm{ if }\, C_1<C_2.
\end{equation*}
Furthermore, $C_2<C_1$ if $i=\nu$ \textup(and so
$\Pi_{\{2,\ldots\nu\}}(\alpha\tens\xi)=0$\textup).
\item
For $I=\{2,\ldots 2\nu-1\}$, we have
\begin{equation*}
k_I^2 = \frac{\prod_{k=2}^{\nu}(\widetilde w_1+\widetilde w_k)}
{\prod_{k=\nu+1}^{2\nu-1}(\widetilde w_1-\widetilde w_k)}.
\end{equation*}
Equality holds iff $\nabla\xi=\Pi_1(\alpha\tens\xi)$ for a $1$-form
$\alpha$ with $\Pi_{\{2,...\nu\}}(\alpha\tens\xi)=0$.
\item
For $\widehat I=\{i\}$ with $i\in\{2,\ldots 2\nu-1\}$, we have
\begin{equation*}
k_I^2 = \frac{\widetilde w_i+\widetilde w_{2\nu+1-i}}
{\widetilde w_i-\widetilde w_1}
\prod_{\substack{j\in J_i^{\max}\\ j\neq 2\nu+1-i}}
\frac{\widetilde w_i+\widetilde w_j}{\widetilde w_i -\widetilde w_{2\nu+1-j}}.
\end{equation*}
Equality holds iff $\nabla\xi=\Pi_i(\alpha\tens\xi)$ for a $1$-form
$\alpha$ with $\Pi_{J_i^{\max}}(\alpha\tens\xi)=0$. Here
$J_i^{\max}=\{i+1,\ldots\nu,2\nu+1-i,\ldots 2\nu-1\}$.
\item
For $I=\{2,\ldots2\nu-2\}$, we have
\begin{equation*}
k_I^2 = \frac{\prod_{k=2}^{\nu}(\widetilde w_{2\nu-1}+\widetilde w_k)}
{(\widetilde w_{2\nu-1}-\widetilde w_1)\prod_{k=\nu+1}^{2\nu-2}
(\widetilde w_{2\nu-1}-\widetilde w_k)}
+\frac{\prod_{k=2}^{\nu}(\widetilde w_1+\widetilde w_k)}
{\prod_{k=\nu+1}^{2\nu-1}(\widetilde w_1-\widetilde w_k)}.
\end{equation*}
Equality holds iff
$\nabla\xi=\Pi_{\{1,2\nu-1\}}(\alpha\tens\xi)$ for a $1$-form $\alpha$
with $\Pi_{\{2,...\nu\}}(\alpha\tens\xi)=0$.
\textup(This is not a refined inequality when $N=3$.\textup)
\item
For $\widehat I=\{i,2\nu-i\}$ with $i\in\{2,\ldots\nu-1\}$, we have
\begin{equation*}
k_I^2 = \frac{\widetilde w_i+\widetilde w_{2\nu+1-i}} {\widetilde
w_i-\widetilde w_1} \prod_{\substack{j\in J_i^{\max}\\ j\neq 2\nu+1-i}}
\frac{\widetilde w_i+\widetilde w_j}{\widetilde w_i -\widetilde w_{2\nu+1-j}}
+\frac{\widetilde w_{i+1}+\widetilde w_{2\nu-i}} {\widetilde
w_{2\nu-i}-\widetilde w_1} \prod_{\substack{j\in J_i^{\max}\\ j\neq i+1}}
\frac{\widetilde w_{2\nu-i}+\widetilde w_j}{\widetilde w_{2\nu-i}-\widetilde
w_{2\nu+1-j}}.
\end{equation*}
Equality holds iff
$\nabla\xi=\Pi_{\{i,2\nu-i\}}(\alpha\tens\xi)$ for a $1$-form $\alpha$
with $\Pi_{J_i^{\max}}(\alpha\tens\xi)=0$. Here
$J_i^{\max}=\{i+1,\ldots\nu,2\nu+1-i,\ldots 2\nu-1\}$.
\end{enumerate}
Replacing $\nu$ by $\nu+1$ gives analogous results for $N=2\nu+1$, but note
that equality cases with $\Pi_{\nu+1}(\alpha\tens v)=0$ will not be attained
if $\lambda$ is properly half-integral.
\end{thm}

We now give more detailed formulas when $N=3$, which is the most common case
arising in practice: the representation $\tau\tens\lambda$
splits into $N=3$ components when:
\begin{enumerate}
\item $V = \odot^k\Lambda^p$ ($k$ a positive integer) and $0<p \leq m-1$
($p=m-1$ in even dimension belongs to this case only by virtue of our
convention on distinctness of conformal weights). Then $\lambda = (k,\ldots
k,0,\ldots0)$ where $k$ is repeated $p$ times and $w_1 = k > w_2 = -p > w_3 =
p-k+1-n$.
\item in odd dimensions, $V=\odot^{k}\Lambda^m$ ($k$ a positive integer) or
$V=\odot^{k-\frac12}\Lambda^m\odot\Delta$ ($k>1/2$ and half-integral), where
$\Delta$ is the spin representation. This corresponds in both cases to
$\lambda = (k,\ldots k)$ and $w_1 = k > w_2 = -\frac{n-1}2 > w_3 =
-k-\frac{n-1}2$.
\end{enumerate}
Note that $P_1$ and $P_2+P_3$ are elliptic, whereas $P_2$ and $P_3$ are
non-elliptic, unless $\nu=1$ and $\lambda$ is properly half-integral, when
$P_2$ is elliptic, but the results above do not cover this case.

\begin{thm} \label{thkatoN3}
If $\xi$ is a nonvanishing section in the kernel of one of the elliptic
operators $P_1$, $P_2+P_3$, $P_1+P_3$ or $P_1+P_2$, we have a refined Kato
inequality $\bigl|d|\xi|\bigr| \leq k_I\, |\nabla\xi|$ with $k_I$ given as
follows.
\begin{enumerate}
\item For $P_1$ or $P_1+P_3$,
$$k_{\{1\}}^2=k_{\{1,3\}}^2=
\frac{w_1}{w_1 - w_2} = \frac{k}{k + p}$$ and equality holds iff
$\nabla \xi = \Pi_2 (\alpha \tens \xi)$ for a $1$-form $\alpha$ such that
$\Pi_3(\alpha \tens \xi) = 0$.
\item For $P_2+P_3$,
$$k_{\{2,3\}}^2=\frac{-w_3}{w_1 - w_3} = \frac{k+n-p-1}{2k+n-p-1}$$
and equality holds iff $\nabla \xi = \Pi_1(\alpha\tens\xi)$ for a
$1$-form $\alpha$ such that $\Pi_2(\alpha\tens\xi)=0$.
\item For $P_1+P_2$,
$$k_{\{1,2\}}^2= \frac{w_1}{w_1 - w_3}= \frac{k}{2k + n - p -1}$$
and equality holds iff $\nabla \xi = \Pi_3 (\alpha \tens \xi)$ for a $1$-form
$\alpha$ such that $\Pi_2(\alpha \tens \xi) = 0$.
\end{enumerate}
\end{thm}
When $\lambda$ is properly half-integral, only the first constant is sharp and
we do not get a nontrivial constant for $P_2$. Since this case sometimes
arises in practice (e.g., the Rarita-Schwinger operator), we note briefly how
the Kato constant can be found. Since $\widetilde w_{\nu+1}=0$, the projection
$\Pi_{\nu+1}=\Pi_2$ is a equal to $\widetilde A_{N-1}=\widetilde A_2$ divided
by $\widetilde w_1\widetilde w_3<0$. Hence we need to obtain a better upper
bound on $\ip{\widetilde A_2(\alpha\tens v),\alpha\tens v}=
\ip{(B^2-\widetilde w_1\widetilde w_3)(\alpha\tens v),\alpha\tens v}$.  Now
for fixed $\alpha\neq0$, say $\alpha=e_n$, we can break this up under
$\lie{so}(n-1)$ and use the fact, easily verified, that $B^2$ is the
difference between the Casimir number of $\lambda$ and the Casimir operator of
$\lie{so}(n-1)$. Applying the branching rule, we see that the eigenvalues of
$(\widetilde A_2)_{e_n\tens e_n}$ are $-(k-\ell)^2$ for $\ell\in\N$ with
$k-\ell\geq0$. Hence if $k$ is half-integral, $\ip{\widetilde A_2(\alpha\tens
v),\alpha\tens v}\leq -1/4$. This gives:
\begin{align*}
k_{\{2\}}^2&=1-\frac1{2k(2k+n-1)}\\
k_{\{2,3\}}&=\frac{(2k+n-1)^2-1}{(2k+n-1)(4k+n-1)}\qquad
k_{\{1,2\}}=\frac{k^2-1}{k(4k+n-1)}
\end{align*}
The analogues of these sharper results for larger $N=2\nu+1$, can be derived
from Branson's minimization formula~\cite{branson-kato}. In particular
he gives the formula for $k_{\{\nu+1\}}$ explicitly there.

Most ``uncomplicated'' tensor bundles, such as vectors, forms, symmetric
traceless tensors and algebraic Weyl tensors, have $N=3$ (except in
low dimensions, where $N$ might be $2$). 
\begin{enumerate}
\item For $\Lambda^1$, the constants are $\frac12$ (conformal or Killing
vector fields), $\frac{n - 1}{n}$ (harmonic $1$-forms) and $\frac1{n}$ (closed
$1$-forms dual to a conformal vector field). The last of these is trivial,
since the only non-vanishing component of $\nabla\xi$ in this case is $\frac1n
\divg\xi\,\iden$.
\item For $\Lambda^2$, the constants are
$\frac13$, $\frac{n-2}{n-1}$ and
$\frac1{n-1}$.
The second of these is the constant for harmonic $2$-forms.
\item For $S^2_0$, the constants are
$\frac23$, $\frac{n}{n+2}$ and $\frac2{n+2}$.
The second of these is the constant appearing in the work of R.~Schoen,
L.~Simon and S.~T.~Yau \cite{ssy}.
\item For $\Lambda^2\odot\Lambda^2$, the constants are
$\frac12$, $\frac{n-1}{n+1}$ and $\frac2{n+1}$.
The second of these is the constant for the second Bianchi identity
appearing in the work of S.~Bando, A.~Kasue and H.~Nakajima~\cite{bkn}.
\end{enumerate}

\section{Refined Kato inequalities with $N$ even}

When $N=2\nu$ is even, we have to minimize or maximize over $J\in\NE$, a sum
of a subset of the following terms:
\begin{align*}
\frac{(\widetilde w_i-\frac12)(\widetilde w_i+\widetilde w_{N+2-i})}
{(\widetilde w_i-\widetilde w_1)(\widetilde w_i-\widetilde w_{\nu+1})}
\prod_{\substack{j\in J\\ j\neq N+2-i}}
\frac{\widetilde w_i+\widetilde w_j}{\widetilde w_i -\widetilde w_{N+2-j}}
&&\qquad\textrm{for}\quad i\in \widehat J\setminus\{1,\nu+1\}\\
\frac{\widetilde w_1-\frac12}{\widetilde w_1-\widetilde w_{\nu+1}}
\prod_{j\in J}\frac{\widetilde w_1+\widetilde w_j}
{\widetilde w_1-\widetilde w_{N+2-j}}&&\\
\frac{\widetilde w_{\nu+1}-\frac12}{\widetilde w_{\nu+1}-\widetilde w_1}
\prod_{j\in J}\frac{\widetilde w_{\nu+1}+\widetilde w_j}
{\widetilde w_{\nu+1}-\widetilde w_{N+2-j}}&&
\end{align*}
Using Proposition~\ref{useful}, the first expression is minimized (subject to
$J\not\ni i$) by $J_i^{\min}=\{2,\ldots i-1,N+2-\nu,\ldots N+2-i\}$ and is
maximized by $J_i^{\max}=\{i+1,\ldots\nu,N+2-i,\ldots N\}$. The second
expression is minimized by $J_1^{\min}=\{N+2-\nu,\ldots N\}$ and maximized by
$J_1^{\max}=\{2,\ldots\nu)$, while the third expression is minimized by
$J_{\nu+1}^{\min}=\{2,\ldots\nu)$ and maximized by
$J_{\nu+1}^{\max}=\{N+2-\nu,\ldots N\}$.

We now proceed as in the odd dimensional case, except that the
analogue of Lemma~\ref{ratio} is no longer useful, due to the additional
$\widetilde w-\frac12$ factors. The results are summarized below.
\begin{thm} Let $E$ be associated to a representation $\lambda$ with $N=2\nu$ 
and let $P_I$ an elliptic operator on sections of $E$ associated to
a subset $I$ of $\{1,\ldots N\}$. Then in the following cases, a refined Kato
inequality of the type $\bigl|d|\xi|\bigr| \leq k_I |\nabla\xi|$ holds outside
the zero set of $\xi$ for $\xi$ in the kernel of $P_I$.
\begin{enumerate}
\item
For $\{1\} \subseteq I \subseteq \{1,\nu+2,\ldots 2\nu\}$, we have
\begin{equation*}
k_I^2 = 1-\frac{(\widetilde w_1-\frac12)
\prod_{k=\nu+2}^{2\nu}(\widetilde w_1+\widetilde w_k)}
{\prod_{k=2}^{\nu+1}(\widetilde w_1-\widetilde w_k)}.
\end{equation*}
Equality holds iff $\nabla\xi=\Pi_{\{2,...\nu+1\}}(\alpha\tens\xi)$ for a
$1$-form $\alpha$ with $\Pi_{\{\nu+2,...2\nu\}}(\alpha\tens\xi) =0$.
\item
For $\{\nu+1\}\subseteq I\subseteq\{2,\ldots\nu,\nu+1\}$, we have
\begin{equation*}
k_I^2 = 1-\frac{(\widetilde w_{\nu+1}-\frac12) \prod_{k=2}^{\nu}(\widetilde
w_{\nu+1}+\widetilde w_k)} {(\widetilde w_{\nu+1}-\widetilde
w_1)\prod_{k=\nu+2}^{2\nu}(\widetilde w_{\nu+1}-\widetilde w_k)}.
\end{equation*}
Equality holds iff $\nabla\xi=\Pi_{\{1,\nu+2...2\nu\}}(\alpha\tens\xi)$ for a
$1$-form $\alpha$ with $\Pi_{\{2,...\nu\}}(\alpha\tens\xi)=0$.
\item
For $\{i,2\nu+2-i\} \subseteq I \subseteq \{i,2\nu+2-i\}\cup J_0$, with $i\in
\{2,\ldots,\nu\}$ and $J_0=\{j:2\leq j<i\}\cup\{2\nu+2-j:i<j\leq\nu\}$, we
have
\begin{equation*}
k_I^2 = 1 - \min (C_1,C_2)
\end{equation*}
where 
\begin{equation*}\begin{split}
C_1 &= \frac{(\widetilde w_i+\widetilde w_{2\nu+2-i})(\widetilde w_i
-\frac12)}{(\widetilde w_i-\widetilde w_{\nu+1})(\widetilde w_i
-\widetilde w_1)} \ \prod_{k\in J_0}\frac{\widetilde w_i+
\widetilde w_k}{\widetilde w_i-\widetilde w_{2\nu+2-k}}, \\
C_2 &= \frac{(\widetilde w_i+\widetilde w_{2\nu+2-i})(\widetilde w_{2\nu+2-i}
-\frac12)}{(\widetilde w_{2\nu+2-i}-
\widetilde w_{\nu+1})(\widetilde w_{2\nu+2-i}-\widetilde w_1)}
\prod_{k\in J_0}\frac{\widetilde w_{2\nu+2-i}+
\widetilde w_k}{\widetilde w_{2\nu+2-i}-\widetilde w_{2\nu+2-k}}.
\end{split}\end{equation*}
Equality holds iff $\nabla\xi=
\Pi_{\widehat J^0\setminus\{i,2\nu+2-i\}}(\alpha\tens\xi)$ for a $1$-form 
$\alpha$ with
\begin{equation*}
\Pi_{\{i\}\cup J_0}(\alpha\tens\xi)=0\,\textrm{ if }\, C_2<C_1\quad
\textrm{or}\quad
\Pi_{\{2\nu+2-i\}\cup J_0}(\alpha\tens\xi)=0\,\textrm{ if }\, C_1<C_2.
\end{equation*}
\item
For $I=\{2,\ldots2\nu\}$, we have
\begin{equation*}
k_I^2 = \frac{(\widetilde w_1-\frac12)
\prod_{k=2}^{\nu}(\widetilde w_1+\widetilde w_k)}
{\prod_{k=\nu+1}^{2\nu}(\widetilde w_1-\widetilde w_k)}.
\end{equation*}
Equality holds iff $\nabla\xi=\Pi_1(\alpha\tens\xi)$ for a $1$-form $\alpha$
with $\Pi_{\{2,...\nu\}}(\alpha\tens\xi)=0$.
\item
For $I=\{1,\ldots\nu,\nu+2,\ldots2\nu\}$, we have
\begin{equation*}
k_I^2 = \frac{(\widetilde w_{\nu+1}-\frac12)
\prod_{k=\nu+2}^{2\nu}(\widetilde w_{\nu+1}+\widetilde w_k)}
{\prod_{k=1}^{\nu}(\widetilde w_{\nu+1}-\widetilde w_k)}.
\end{equation*}
Equality holds iff $\nabla\xi=\Pi_{\nu+1}(\alpha\tens\xi)$ for a $1$-form
$\alpha$ with $\Pi_{\{\nu+2,...2\nu\}}(\alpha\tens\xi)=0$.
\item
For $\widehat I=\{i\}$ with $i\in\{2,\ldots\nu,\nu+2,\ldots 2\nu\}$, we have
\begin{equation*}
k_I^2 = \frac{(\widetilde w_i-\frac12)(\widetilde w_i+\widetilde w_{2\nu+2-i})}
{(\widetilde w_i-\widetilde w_1)(\widetilde w_i-\widetilde w_{\nu+1})}
\prod_{\substack{j\in J_i^{\max}\\ j\neq 2\nu+2-i}}
\frac{\widetilde w_i+\widetilde w_j}{\widetilde w_i-\widetilde w_{2\nu+2-j}}.
\end{equation*}
Equality holds iff $\nabla\xi=\Pi_i(\alpha\tens\xi)$ for a $1$-form
$\alpha$ with $\Pi_{J_i^{\max}}(\alpha\tens\xi)=0$. Here
$J_i^{\max}=\{i+1,\ldots\nu,2\nu+2-i,\ldots 2\nu\}$.
\item
For $I=\{2,\ldots2\nu-1\}$, we have
\begin{equation*}
k_I^2 = \frac{(\widetilde w_1-\frac12) \prod_{k=2}^{\nu}(\widetilde
w_1+\widetilde w_k)} {\prod_{k=\nu+1}^{2\nu}(\widetilde w_1-\widetilde w_k)}+
\frac{(\widetilde w_{2\nu}-\frac12) \prod_{k=2}^{\nu}(\widetilde
w_{2\nu}+\widetilde w_k)} {(\widetilde w_{2\nu}-\widetilde
w_1)\prod_{k=\nu+1}^{2\nu-1}(\widetilde w_{2\nu}-\widetilde w_k)}.
\end{equation*}
Equality holds iff $\nabla\xi=\Pi_{\{1,2\nu\}}(\alpha\tens\xi)$ for a $1$-form
$\alpha$ with
$\Pi_{\{2,...\nu\}}(\alpha\tens\xi)=0$.
\item
For $I=\{1,\ldots\nu-1,\nu+2,\ldots2\nu\}$ we have
\begin{equation*}
k_I^2 = \frac{(\widetilde w_\nu-\frac12) \prod_{k=\nu+2}^{2\nu}(\widetilde
w_\nu+\widetilde w_k)} {(\widetilde w_\nu-\widetilde
w_{\nu+1})\prod_{k=1}^{\nu-1}(\widetilde w_\nu-\widetilde w_k)}
+\frac{(\widetilde w_{\nu+1}-\frac12) \prod_{k=\nu+2}^{2\nu}(\widetilde
w_{\nu+1}+\widetilde w_k)} {\prod_{k=1}^{\nu}(\widetilde w_{\nu+1}-\widetilde
w_k)}.
\end{equation*}
Equality holds iff $\nabla\xi=\Pi_{\{\nu,\nu+1\}}(\alpha\tens\xi)$ for a
$1$-form $\alpha$ with $\Pi_{\{\nu+2,...2\nu\}}(\alpha\tens\xi)=0$.
\item
For $\widehat I=\{i,2\nu+1-i\}$ with $i\in\{2,\ldots\nu-1\}$, we have
\begin{multline*}
k_I^2 = \frac{(\widetilde w_i+\widetilde w_{2\nu+2-i})(\widetilde
w_i-\frac12)}{(\widetilde w_i-\widetilde w_1)(\widetilde w_i-\widetilde
w_{\nu+1})} \prod_{\substack{j\in J_i^{\max}\\ j\neq 2\nu+2-i}}
\frac{\widetilde w_i+\widetilde w_j}{\widetilde w_i -\widetilde
w_{2\nu+2-j}}\\
+\frac{(\widetilde w_{i+1}+\widetilde w_{2\nu+1-i})(\widetilde w_i-\frac12)}
{(\widetilde w_{2\nu+1-i}-\widetilde w_1)(\widetilde w_{2\nu+1-i}-\widetilde
w_{\nu+1})} \prod_{\substack{j\in J_i^{\max}\\ j\neq i+1}} \frac{\widetilde
w_{2\nu+1-i}+\widetilde w_j}{\widetilde w_{2\nu+1-i}-\widetilde w_{2\nu+2-j}}.
\end{multline*}
Equality holds iff $\nabla\xi=\Pi_{\{i,2\nu+1-i\}}(\alpha\tens\xi)$ for a
$1$-form $\alpha$ with $\Pi_{J_i^{\max}}(\alpha\tens\xi)=0$. Here
$J_i^{\max}=\{i+1,\ldots\nu,2\nu+2-i,\ldots 2\nu\}$.
%
%
\end{enumerate}
\end{thm}
We now give more detailed formulas when $N=4$, which is the generic
case in four dimensional differential geometry: the representation
$\tau\tens\lambda$ splits into $N=4$ components whenever
\begin{enumerate}
\item if $n=2m$ is even and
$V=\odot^{\ell}\Lambda_{\pm}\odot^{k-\ell}\Lambda^p$ or
$V=\odot^{\ell-\frac12}\Lambda_{\pm} \odot^{k-\ell}\Lambda^p\odot\Delta_{\pm}$
where $k>\ell >0$ are (simultaneously) integers or half-integers, $p<m$ are
integers, $\Lambda^m_{\pm}$ stand for selfdual or antiselfdual $m$-forms and
$\Delta_{\pm}$ for positive or negative spin representations.  The associated
weights are $\lambda = (k,\ldots k,\ell,\ldots\ell,\pm\ell)$, with $k$
repeated $p$ times.  One gets $ w_1 = k > w_2 = \ell - p > w_3 = 1 -\frac{n}2
- \ell > w_4 = - k+p+1-n$.
\item if $n=2m+1$ is odd, $V=\odot^{k-\frac12} \Lambda^p \odot \Delta$ with
$p<m$ integer and $k\geq\frac12$ and half-integer, so that $\lambda =
(k,\ldots k, \frac12, \ldots \frac12)$. Conformal weights are a specialization
of the previous formula with $\ell=\frac12$: $w_1 = k > w_2 = \frac12 - p >
w_3 = \frac{(1 - n)}2 > w_4 = - k + p + 1 - n$.
\end{enumerate}

Note that that $P_1$, $P_3$ and $P_2+P_4$ are elliptic, whereas
$P_2$ and $P_4$ are non-elliptic.

We give in the following theorem the Kato constants for the kernels of the
minimal elliptic operators.
\begin{thm}\label{thkatoN4}
If $\xi$ is a nonvanishing section in the kernel of one of the elliptic
operators $P_1$, $P_3$ or $P_2+P_4$, we have a refined Kato inequality
$\bigl|d|\xi|\bigr| \leq k_I\, |\nabla\xi|$ with $k_I$ given as follows.
\begin{enumerate}
\item For $P_1$,
\begin{equation*} \begin{split}
k_{\{1\}}^2
& = 1 - \frac{(w_1+\frac{n-2}2)(w_1+w_4+n-1)}
{(w_1 - w_2)(w_1 - w_3)}\\ 
& = \frac{(k + \frac{n - 2}2)(k-\ell) + \ell (k - \ell + p)}
{(k- \ell + p) (k + \ell + \frac{n - 2}2)} 
\end{split}\end{equation*}
and equality holds iff
$\nabla \xi = (\Pi_2+\Pi_3)(\alpha \tens \xi)$ for a
$1$-form $\alpha$ with $\Pi _4 (\alpha \tens \xi) = 0$.
\item For $P_3$,
\begin{equation*} \begin{split} 
k_{\{3\}}^2 &= 1 - \frac{(w_3+\frac{n-2}2)(w_3+w_2+n-1)}{(w_3-w_4)(w_3-w_1)}\\
&= 1 - \frac{\ell(\frac{n}2-p)}{(k-\ell+\frac{n}2-p)(k+\ell+\frac{n-2}2)}
\end{split} \end{equation*}
and equality holds iff $\nabla \xi = (\Pi_1+\Pi_4)(\alpha \tens \xi)$ for a
$1$-form $\alpha$ with $\Pi_2 (\alpha \tens \xi) = 0$.
\item For $P_2+P_4$,
\begin{equation*} \begin{split}
k_{\{2,4\}}^2 & = 1-\min\left\{\frac{(w_4+\frac{n-2}2)(w_2+w_4+n-1)}
{(w_4 - w_1) (w_4 - w_3)},\frac{(w_2+\frac{n-2}2)(w_2+w_4+n-1)}
{(w_2 - w_1) (w_2 - w_3)}\right\}\\
& = 1-\min\left\{\frac{(k+\frac n2-p)(k-l)}
{(2k+n-p-1) (k-l+\frac n2-p)},\frac{(l+\frac n2-p-1)(k-l)}
{(k-l+p)(2l+\frac n2-p-1)}
\right\}
\end{split}\end{equation*}
and equality holds iff $\nabla\xi = (\Pi_1+\Pi_3)(\alpha \tens \xi)$ for a
$1$-form $\alpha$ with $\Pi_2 (\alpha\tens\xi)=0$ or
$\Pi_4(\alpha\tens\xi)=0$ depending on which term is the minimum.
\end{enumerate}
\end{thm}

\section*{Appendix: Explicit constants for dimension $3$ and $4$}

{\flushleft\bf Dimension $3$}. Irreducible representations of $\lie{so}(3)$
are symmetric powers, denoted $\Delta^r$, of the spin representation $\Delta$
(if $r$ is even, $\Delta^r$ has a canonical real structure and we denote from
now on by $\Delta^r$ its real part). The Clebsch-Gordan formulas show that we
are in the case $N=2$ if $r=1$ and $N=3$ if $r\geq 2$. In the former case, the
elliptic operators are the (Penrose) twistor operator $P_1$ and the Dirac
operator $P_2$ corresponding to projections on the first and second part of
\begin{equation*}
\R^3 \tens \Delta = \Delta^2  \tens \Delta = \Delta^3 \oplus \Delta.
\end{equation*}
In the latter case, the elliptic operators are the twistor operator $P_1$ and
Dirac-type operator $P_2+P_3$ corresponding to projections on the first or
second-and-third part of
\begin{equation*}
\R^3 \tens \Delta^r = \Delta^2 \tens \Delta^r  
= \Delta^{r+2} \oplus \Delta^{r} \oplus \Delta^{r-2}, \quad r\geq 2.
\end{equation*}
If $r\geq3$ and is odd, then $P_2$ is elliptic on its own: it is the
Rarita-Schwinger operator when $r=3$ and so we denote it by R-S in general.
 
The following table sums up our formulae in three dimensions.

\medskip

\begin{small}
\begin{center}
\begin{tabular}{|c|c|c|}
\hline
operator & conditions & refined constant \\ 
\hline
Twistor & all $r$ & $\sqrt{\frac{r}{r+2}}$ \\
\hline
Dirac & $r=1$ & $\sqrt{\frac23}$ \\
\hline
Dirac-type & $r\geq 2$ & $\sqrt{\frac{r+2}{2(r+1)}}$ \\
\hline
R-S ($r$ odd) & $r\geq 3$ & $\sqrt{1-\frac{1}{r(r+2)}}$ \\
\hline
\end{tabular}
\end{center}
\end{small}

\medskip

{\flushleft\bf Minimal elliptic operators in dimension $4$}.  Irreducible
representations of $\mathfrak{so}(4)$ are tensor products of symmetric powers,
denoted $V^{r,s}= \Delta_{+}^r\tens\Delta_{-}^r$, of the positive and negative
half-spin representation $\Delta_{\pm}$ (if $r+s$ is even, $V^{r,s}$ has a
canonical real structure and, as above, $V^{r,s}$ will denote its real part).
Assuming $r \geq s$, the Clebsch-Gordan formulas yield, for $r \geq s > 0$,
\begin{equation*} 
\R^4 \tens V^{r,s}
= V^{r+1,s+1} \oplus V^{r+1,s-1} \oplus V^{r-1,s+1} \oplus V^{r-1,s-1}
\end{equation*}
so that we are in the case $N=4$ if $r>s>0$ and the case $N=3$ if $r = s > 0$
(the middle components have equal conformal weights here). If $r > s = 0$ then
\begin{equation*}
\R^4 \tens V^{r,0} = V^{r+1,1} \oplus V^{r-1,1},
\end{equation*}
and we are in the case $N=2$.

Hence we have (at most) three minimal elliptic operators.

\begin{enumerate}

\item The twistor operator, given by the projection on the first factor in
every case.

\smallskip

\item[(ii-a)] The operator given by from the projection onto $V^{r-1,s+1}$.
It is the operator $P_3$ when $N=4$ (\ie, if $r>s>0$) or $P_2$ when $N=2$
(\ie, if $s$~vanishes). It defines the spin $\frac r2$ field equation in this
last case and we shall call it a ``spin $\frac{r+s}2$ field'' in general.

\smallskip

\item[(ii-b)] The operator in (ii-a) is not elliptic if $N=3$ (\ie, if
$r=s>0$). We shall replace it by the one given by the projection onto
$\left(V^{s+1,s-1} \oplus V^{s-1,s+1}\right) \oplus V^{s-1,s-1}$. The usual
Hodge-de Rham belongs to this case, so that it seems reasonable to call it
a Dirac-type operator.

\smallskip

\item[(iii)] The operator given by the projection onto $V^{r+1,s-1}\oplus
V^{r-1,s-1}$ is the elliptic operator $P_2+P_4$ if $N=4$ (\ie, if $r>s>0$).
We shall again call it a Dirac-type operator.

\end{enumerate}
 
The following table sums up our formulae in four dimensions.

\medskip

\begin{small}
\begin{center}
\begin{tabular}{|c|c|c|c|c|}
\hline
operator & conditions & refined constant & s=0 & r=s \\ 
\hline
Twistor & $r\geq s\geq 0$ & $\sqrt{\frac{2rs+r+s}{2(r+1)(s+1)}}$
& $\sqrt{\frac r{2(r+1)}}$ & $\sqrt{\frac s{s+1}}$\\
\hline
Spin $\frac{r+s}2$ field& $r>s\geq0$ & $\sqrt{\frac{2rs+r+3s+2}{2(r+1)(s+1)}}$
& $\sqrt{\frac {r+2}{2(r+1)}}$ & -\\
\hline
Dirac-type & $r\geq s>0$ & $\sqrt{\frac{s+2}{2(s+1)}}$
& - & $\sqrt{\frac {s+2}{2(s+1)}}$\\
\hline 
\end{tabular}
\end{center}
\end{small}

\smallskip

As an example, we can obtain from the table the value found by M.~Gursky and
C.~LeBrun in \cite{gursky-lebrun-einstein} for a co-closed positive half Weyl
tensor (outside its zero set):
\begin{equation}
\bigl|d|W^+|\bigr| \leq \sqrt{\frac35} \, |\nabla W^+|,
\end{equation}
and notice that equality occurs if and only if $\nabla W^+ = \Pi_2
(\alpha\tens W^+)$.

\vskip.5cm

\end{document}